\documentclass[aop,citesort,MSNbibl,dvips]{arximspdf}

\doi{10.1214/10-AOP603}
\volume{39}
\issue{6}
\pubyear{2011}
\firstpage{2441}
\lastpage{2473}

\makeatletter

\newcommand{\underlineIk}{\underline{i}_{\hspace*{0.8pt}k}}
\newcommand{\underlineJk}{\underline{j}_{\hspace*{0.8pt}k}}

\newtheorem{Theorem}{Theorem}
\newtheorem{Lemma}{Lemma}
\newtheorem{Corollary}{Corollary}
\newtheorem{Proposition}{Proposition}

\newproclaim{Definition}{Definition}
\newproclaim{Assumption}{Assumption}

\newcommand{\convd}{\stackrel{\mathcal{D}}{\longrightarrow}}
\newcommand{\convP}{\stackrel{\mathcal{P}}{\longrightarrow}}

\newcommand{\Ykm}{Y_{k}^{(m)}}
\newcommand{\ld}{\mathcal{M}}
\newcommand{\corr}{\operatorname{Corr}}

\newcommand{\Var}{\operatorname{Var}}
\newcommand{\Cov}{\operatorname{Cov}}

\makeatother

\begin{document}
\begin{frontmatter}

\title{Split invariance principles for stationary processes}
\runtitle{Split invariance principles}

\begin{aug}
\author[A]{\fnms{Istv\'{a}n} \snm{Berkes}\corref{}\thanksref{t1}\ead[label=e1]{berkes@tugraz.at}},
\author[B]{\fnms{Siegfried} \snm{H\"{o}rmann}\thanksref{t2}\ead[label=e2]{shormann@ulb.ac.be}} and
\author[A]{\fnms{Johannes} \snm{Schauer}\thanksref{t3}\ead[label=e3]{johannes@schauer.com}}
\runauthor{I. Berkes, S. H\"{o}rmann and J. Schauer}
\affiliation{Graz University of Technology, Universit\'{e} Libre de
Bruxelles and\break Graz University of Technology}
\address[A]{I. Berkes\\
J. Schauer\\
Graz University of Technology\\
M\"{u}nzgrabenstrasse 11\\
8010 Graz\\
Austria\\
\printead{e1}\\
\hphantom{E-mail: }\printead*{e3}} 
\address[B]{S. H\"{o}rmann\\
Universit\'{e} Libre de Bruxelles\\
Bd du Triomphe\\
1050 Bruxelles\\
Belgium\\
\printead{e2}}
\end{aug}

\thankstext{t1}{Supported by FWF Grant S9603-N23 and OTKA Grants K
67961 and K 81928.}

\thankstext{t2}{Supported by the Banque National de Belgique and
Communaut\'{e} fran\c{c}aise de Belgique---Actions de Recherche
Concert\'{e}es (2010--2015).}

\thankstext{t3}{Supported by FWF Grant S9603-N23.}

\received{\smonth{12} \syear{2009}}
\revised{\smonth{6} \syear{2010}}

%
\begin{abstract}
The results of Koml\'{o}s, Major and Tusn\'{a}dy give optimal Wiener
approximation of partial sums of i.i.d. random variables and provide
an extremely powerful tool in probability and statistical inference.
Recently Wu [\textit{Ann. Probab.} \textbf{35} (2007) 2294--2320]
obtained Wiener approximation of a class of
dependent stationary processes with finite $p$th moments, $2<p\le4$,
with error term $o(n^{1/p}(\log n)^\gamma)$, $\gamma>0$, and Liu and
Lin [\textit{Stochastic Process. Appl.} \textbf{119} (2009) 249--280]
removed the logarithmic factor, reaching the
Koml\'{o}s--Major--Tusn\'{a}dy bound $o(n^{1/p})$. No similar results
exist for $p>4$, and in fact, no existing method for dependent
approximation yields an a.s. rate better than $o(n^{1/4})$. In this
paper we show that allowing a second Wiener component in the
approximation, we can get rates near to $o(n^{1/p})$ for arbitrary
$p>2$. This extends the scope of applications of the results
essentially, as we illustrate it by proving new limit theorems for
increments of stochastic processes and statistical tests for short term
(epidemic) changes in stationary processes. Our method works under a
general weak dependence condition covering wide classes of linear and
nonlinear time series models and classical dynamical systems.
\end{abstract}

%
\begin{keyword}[class=AMS]
\kwd{60F17}
\kwd{60G10}
\kwd{60G17}.
\end{keyword}
\begin{keyword}
\kwd{Stationary processes}
\kwd{strong invariance principle}
\kwd{KMT approximation, dependence}
\kwd{increments of partial sums}.
\end{keyword}

\end{frontmatter}

\section{Introduction}\label{se:intro}

Let $X, X_1, X_2, \ldots$ be i.i.d. random variables with mean~0 and
variance 1, and
let $S_n=\sum_{k\le n} X_k$. Koml\'{o}s, Major and Tusn\'{a}dy \cite
{kmt75,kmt76} showed that
if
$E(e^{t|X|})<\infty$ for some $t>0$ then, after suitably enlarging the
probability space, there exists a Wiener process $\{W(t), t\ge0\}$
such that
%
%
\begin{equation}\label{kmt1}
S_n=W(n)+O(\log n) \qquad\mbox{a.s.}
\end{equation}
Also, if $E|X|^p<\infty$ for some $p>2$, they proved the approximation
%
%
\begin{equation}\label{kmt2}
S_n=W(n)+o(n^{1/p}) \qquad\mbox{a.s.}
\end{equation}
The remainder terms in (\ref{kmt1}) and (\ref{kmt2}) are optimal. In
the case when only $EX=0$, $EX^2=1$ is assumed, Strassen \cite{str64}
obtained 
%
%
\begin{equation} \label{str}
S_n=W(n)+o((n\log\log n)^{1/2}) \qquad\mbox{a.s.}
\end{equation}
Without additional moment assumptions the rate in (\ref{str}) is also
optimal (see Major~\cite{major76}).
Relation (\ref{str}) is a useful invariance principle for the law of
the iterated logarithm; on the other hand, it does not imply the CLT
for $\{X_n\}$. This difficulty was removed by Major \cite{maj} who
showed that under $EX=0$, $EX^2=1$ there exists a Wiener process $W$
and a numerical sequence $\tau_n \sim n$ such that
%
%
\begin{equation}\label{ma1}
S_n=W(\tau_n)+o(n^{1/2}) \qquad\mbox{a.s.}
\end{equation}
Thus allowing a slight perturbation of the approximating Wiener process
one can reach the remainder term $o(n^{1/p})$ also for $p=2$, making
the result applicable for a wide class of CLT-type results.
The case of strong approximation under the moment condition $EX^2
h(|X|)<\infty$ where $h(x)=o(x^\varepsilon)$, $x\to\infty$, for any
$\varepsilon>0$, has been cleared up completely by Einmahl \cite{einmahl}.

The previous results, which settle the strong approximation problem for
i.i.d. random variables with finite variances, provide powerful tools
in probability and statistical inference (see, e.g., Shorack and
Wellner \cite{showe}). Starting with Strassen~\cite{str67}, a wide
literature has dealt with extensions of the above results for weakly
dependent sequences, but the existing results are much weaker than in
the i.i.d. case. Recently, however, Wu \cite{wuap} showed that for a
large class of weakly dependent stationary sequences $\{X_n\}$
satisfying $E|X_1|^p<\infty$, $2<p\le4$, we have the approximation
\[
S_n=W(n)+o(n^{1/p}(\log n)^\gamma) \qquad\mbox{a.s.}
\]
for some $\gamma>0$, and Liu and Lin \cite{liulin}
removed the logarithmic factor in the error term, reaching the optimal
Koml\'{o}s--Major--Tusn\'{a}dy bound. The proofs do not work for $p>4$,
and in fact, no existing method for dependent approximation yields an
a.s. rate better than $o(n^{1/4})$.
On the other hand, many important limit theorems in probability and
statistics involve norming sequences smaller than $n^{1/4}$, making
such results inaccessible by invariance methods. The purpose of the
present paper is to fill this gap and provide a new type of
approximation theorem reaching nearly the Koml\'{o}s--Major--Tusn\'
{a}dy rate for any $p>2$.

As noted above, reaching the error term $o(n^{1/2})$ for i.i.d.
sequences with finite variance requires a perturbation of the
approximating Wiener process~$W$. In the case of dependent processes we
will also need a similar perturbation, and, more essentially, we will
include a second Wiener process in the approximation, whose scaling
factor is smaller than that of $W$, and thus it will not affect the
asymptotic behavior of the main term. Specifically, for a large class
of weakly dependent stationary processes $\{Y_k\}$ with finite $p$th
moments, $2<p<\infty$, we will prove the approximation
%
%
\begin{equation}\label{ketw}
\sum_{k=1}^n Y_k=W_1(s_n^2)+W_2(t_n^2)
+O\bigl(n^{({1+\eta})/{p}}\bigr) \qquad\mbox{a.s.},
\end{equation}
where $\{W_1(t), t\geq0\}$ and $\{W_2(t), t\geq0\}$
are standard Wiener processes, and~$s_n$, $t_n$ are numerical sequences with
\[
s_n^2\sim\sigma^2 n,\qquad t_n^2 \sim c n^\gamma
\]
for some\vspace*{2pt} $0<\gamma<1$, $\sigma^2>0$, $c>0$.
The new element in (\ref{ketw}) is the term $W_2(t_n^2)$ which,
by its smaller scaling, does not disturb the asymptotic properties of~%
$W_1(s_n^2)$. Note that the processes $W_1$, $W_2$ are not independent,
but this will not present any difficulties in applications. (See also
Proposition \ref{pr:corr} in the next section.) The number $\eta$
depends on the weak dependence rate of~$\{Y_k\}$ (introduced below),
and can be made arbitrarily small under suitable rate conditions.

For $p>0$ and a random variable $Y$, let $\|Y\|_p=(E|Y|^p)^{1/p}$.
If $A$ and
$B$ are subsets of $\mathbb{Z}$, we let $d(A,B)=\inf\{|a-b|\dvtx a\in
A,b\in B\}$.
\begin{Definition}\label{d}
Let $\{Y_k, k\in{\mathbb Z}\}$ be a stochastic process, let
$p\ge1$ and let $\delta(m)\to0$. We say that $\{Y_k, k\in{\mathbb
Z}\}$ is \textit{weakly $\mathcal{M}$-dependent in $L^p$ with rate
function $\delta(\cdot)$} if:

\begin{longlist}[(A)]
\item[(A)] For any $k\in{\mathbb Z}$, $m\in{\mathbb N}$ one
can find a random variable $Y_{k}^{(m)}$ with finite $p$th moment such that
\[
\bigl\|Y_k-Y_k^{(m)}\bigr\|_p\le\delta(m).
\]

\item[(B)] For any disjoint intervals $I_1, \ldots, I_r$
($r\in\mathbb{N}$) of integers and any positive integers $m_1,
\ldots, m_r$, the vectors\vspace*{2pt} $\{Y_j^{(m_1)}, j\in I_1\}, \ldots,
\{Y_j^{(m_r)}, j\in I_r\}$ are independent provided $d(I_k,I_l)>\max
\{
m_k,m_l\}$ for $1\leq k<l\leq r$.
\end{longlist}
\end{Definition}
%

We remark that our dependence condition is naturally preserved under
smooth transformations. For example,
if $\{Y_k\}$ is weakly $\mathcal{M}$-dependent in $L^p$ with rate
$\delta(\cdot)$, and $h$ is a Lipschitz $\alpha$
function ($0<\alpha\le1$) with Lipschitz constant $K$, then by the
monotonicity of $\|Y_k-Y_{k}^{(m)}\|_p$ in $p$ we have
%
\[
\bigl\|h(Y_k)-h\bigl(Y_{k}^{(m)}\bigr)\bigr\|_p\leq
K\bigl\|Y_k-Y_{k}^{(m)}\bigr\|_{\alpha p}^{\alpha}\le K\bigl\|Y_k-Y_{k}^{(m)}\bigr\|
_p^{\alpha},
\]
and thus $\{h(Y_k)\}$ is also weakly $\mathcal{M}$-dependent in $L^p$
with rate function~$K\delta(\cdot)^\alpha$.

Note that (B) implies that for any fixed $m$ the
sequence $\{\Ykm, k\in\mathbb{Z}\}$ is an $m$-dependent process.
Hence, sequences satisfying conditions (A) and~(B) are approximable, in
the $L^p$ sense, by $m$-dependent processes of any fixed order $m\geq
1$ with termwise approximation error $\delta(m)$.
In other words, sequences in Definition~\ref{d} are close to $m$-dependent
sequences, the value of $m$ depending on the required closeness,
explaining the terminology.
Since $\|Y_k\|_p\leq\|Y_k^{(m)}\|_p+\|Y_k-Y_k^{(m)}\|_p$, condition
(A) implies
that $E|Y_k|^p$ is finite.
Using $L^p$-distance is convenient for our theorems, but, depending on
the application, other distances can be used in part (A) of
Definition \ref{d}.
For example, defining (as usual) the $L_0$ norm of a random variable
$X$ by
\[
\|X\|_0=\inf\{\varepsilon>0\dvtx P(|X|\ge\varepsilon)<\varepsilon\},
\]
condition (A) could be replaced by
\[
\bigl\|Y_k-Y_k^{(m)}\bigr\|_0\le\delta(m).
\]
Such a definition requires no moment assumptions and turns out to
provide a useful
dependence measure
for studying empirical processes (see \cite{behosch}).

Trivially the previous definition covers $m$-dependent processes for
any fixed $m$
(see also Section \ref{ss:mdep}), but, in contrast to the very
restrictive condition of $m$-dependence, weak $\mathcal{M}$-dependence
holds for a huge class of stationary sequences, including those studied
in Wu \cite{wu05,wuap} and Liu and Lin \cite{liulin}. In the
case when $\{Y_k, k\in{\mathbb Z}\}$ allows a Wiener--Rosenblatt representation
%
%
\begin{equation}\label{eq:Y_k_represent}
Y_k=f(\varepsilon_{k},\varepsilon_{k-1},
\ldots),\qquad k\in{\mathbb Z},
\end{equation}
with an i.i.d. sequence $\{\varepsilon_k, k\in\mathbb{Z}\}$,
weak $\mathcal{M}$-dependence is very close to Wu's
physical dependence condition in \cite{wu05}, except that we allow a
larger freedom in choosing the approximating random variables $\Ykm$,
compared with the choice in \cite{wu05,wuap} via coupling.
(For sufficient criteria for the representation~(\ref
{eq:Y_k_represent}), see Rosenblatt
\cite{rosenblatta,rosenblattb,rosenblatt}.) Note that instead of~(\ref
{eq:Y_k_represent}) we may also assume a~two-sided representation
%
%
\begin{equation}\label{eq:Y_k_represent2}
Y_k=f(\ldots,\varepsilon_{k-1},\varepsilon_k,
\varepsilon_{k+1},\ldots),\qquad k\in{\mathbb Z},
\end{equation}
of $\{Y_k\}$. In case when $\{Y_k, k\in{\mathbb Z}\}$ allows the
representation (\ref{eq:Y_k_represent2}) with mixing~$\{\varepsilon
_k\}
$, Definition \ref{d} is a modified version of NED (see
Section \ref
{ss:ned}), a~weak dependence condition which appeared already in
Ibragimov \cite{ibr} and has been brought forward in Billingsley \cite
{billingsley68} (see also \cite{mcleish75a,mcleish75b}).
Later NED has been successfully used in the econometrics literature to
establish weak dependence of dynamic time series models (see, e.g.,
\cite{pp97}). In Section \ref{se:ex} we will discuss further the
connection between weak $\mathcal{M}$-dependence with known weak
dependence conditions.
We stress that the definition of weak ${\ld}$-de\-pendence does not
assume the representation (\ref{eq:Y_k_represent}) or (\ref
{eq:Y_k_represent2}), although it was motivated by this case. The
reason for using our more general definition is to illuminate the
essential structural condition on $\{Y_k\}$ required for our theorems.
Extensions of our results for ``classical'' mixing conditions, like
$\alpha$, $\beta$, $\rho$ mixing
and their variants will be given in a subsequent paper.

The main results of our paper are formulated in Section \ref{se:main}.
In Section~\ref{se:ex} we give several examples. Applications of the
theorems can be found in
Sections~\ref{se:app} and \ref{se:epi}, while Section \ref
{se:proofs} contains
the proofs of the main theorems.

\section{Main theorems}\label{se:main}
We write $a_n\ll b_n$ if $\overline{\lim}_{n \to\infty
}|a_n/b_n|<\infty$.
\begin{Theorem}\label{th:poly}
Let $p>2$, $\eta>0$ and let $\{Y_k,k\in\mathbb{Z}\}$ be a centered
stationary sequence, weakly $\mathcal M$-dependent in $L^p$ with rate
function
%
%
\begin{equation}\label{eq:ma} \delta(m)\ll m^{-A},
\end{equation}
where
%
%
\begin{equation} \label{Acond}
A>\frac{p-2}{2\eta}\biggl(1-\frac{1+\eta}{p}\biggr)\vee1,\qquad
(1+\eta)/p<1/2.
\end{equation}
Then the series
%
%
\begin{equation}\label{eq:sig}
\sigma^2=\sum_{k\in\mathbb{Z}}EY_0Y_k
\end{equation}
is absolutely convergent, and $\{Y_k,k\in\mathbb{Z}\}$ can be
redefined on a new probability space together with two
Wiener processes $\{W_1(t), t\geq0\}$ and $\{W_2(t), t\geq0\}$
such that
%
%
\begin{equation}\label{eq:strappr}
\sum_{k=1}^n Y_k=W_1(s_n^2)+W_2(t_n^2)
+O\bigl(n^{({1+\eta})/{p}}\bigr) \qquad\mbox{a.s.},
\end{equation}
where $\{s_n\}$ and $\{t_n\}$ are nondecreasing numerical sequences with
%
%
\begin{equation}\label{eq:sn}
s_n^2\sim\sigma^2 n,\qquad t_n^2 \sim c n^\gamma
\end{equation}
for some $0<\gamma<1$, $c>0$.
\end{Theorem}

Note that for any fixed $p>2$ and $0<\eta< (p-2)/2$, condition (\ref
{Acond}) is satisfied if $A$ is
large enough, and thus Theorem \ref{th:poly} provides an a.s.
invariance principle with remainder term close
to the optimal remainder term $o(n^{1/p})$ in the Koml\'
{o}s--Major--Tusn\'{a}dy approximation.

It is natural to ask if $W_1(s_n^2)$ in (\ref{eq:strappr}) can be
replaced by $W_1(\sigma^2 n)$, a fact that would simplify applications.
The proof of the theorem yields an $s_n$ with $s_n^2=\sigma^2
n+O(n^{1-\epsilon})$ for some $0<\epsilon<1$,
but for $A$ barely exceeding the lower bound in (\ref{Acond}), the
explicit value of $\epsilon$ is very small.
Thus replacing\vspace*{1pt} $W_1(s_n^2)$ by $W_1(\sigma^2 n)$ introduces an
additional error term that ruins the error term $O(n^{(1+\eta)/p})$ in
(\ref{eq:strappr}). The situation is similar to the Wiener
approximation of partial sums of i.i.d. random variables\vadjust{\goodbreak}
with mean 0 and variance 1 when we have (\ref{ma1}) with a numerical
sequence $\tau_n \sim n$, but in general (\ref{ma1}) does not hold with
$\tau_n=n$. (See Major \cite{major76,maj}.) Note, however, that
in our case the large difference between $s_n^2$ and $\sigma^2 n$ is a
consequence of the method, and we do not claim that another
construction cannot yield the approximation~(\ref{eq:strappr}) with
$s_n^2=\sigma^2 n$. However, the presence of $s_n^2$ in~(\ref
{eq:strappr}) does not limit the applicability of our strong invariance
principle: $s_n^2$ and $t_n^2$ are explicitly calculable nonrandom
numbers and as we will see, applying limit theorems for $W_1(s_n^2)$ is
as easy as for $W_1(\sigma^2 n)$.

As the proof of Theorem~\ref{th:poly} will show, the sequences $\{s_n\}
$ and $\{t_n\}$ in~(\ref{eq:strappr}) have a complementary character.
More precisely, there is a partition $\mathbb N=G_1\cup G_2$ (provided
by the long and short blocks in a traditional blocking argument) and a
representation
\[
s_n^2=\sum_{k=1}^n \sigma_k^2,\qquad t_n^2=\sum_{k=1}^n \tau_k^2\qquad
(n=1, 2, \ldots)
\]
such that $\sigma_k^2$ converges to $\sigma^2$ on $G_1$ and equals 0 on
$G_2$, and $\tau_k^2$ converges to $\sigma^2$ on $G_2$ and equals 0 on
$G_1$. In particular,
%
%
\begin{equation}\label{stlimsup}
\mathop{\overline{\lim}}_{n\to\infty}
(s_{n+1}^2-s_n^2)=\mathop{\overline{\lim}}
_{n\to\infty} (t_{n+1}^2-t_n^2)=\sigma^2,
\end{equation}
and both liminf's are equal to 0.

The numerical value of $\gamma$ in (\ref{eq:sn}) plays no role in the
applications in this paper, but for later applications we note that
if
\[
A>\frac{p-2}{2\eta
(1-\varepsilon_0)^2}\biggl(1-\frac{1+\eta}{p}\biggr)\vee1
\]
for some $0<\varepsilon_0<1$, then we can choose
%
%
\begin{equation}\label{num}
\gamma=1-\varepsilon_0 \frac{2\eta(1-\varepsilon_0)}{p-2(1+\eta
\varepsilon_0)}.
\end{equation}


As we already mentioned in the \hyperref[se:intro]{Introduction}, the processes $W_1$ and
$W_2$ are not independent. While for our applications this is not important,
the following proposition might be useful for possible further applications.

\begin{Proposition}\label{pr:corr}
Under the assumptions of Theorem \ref{th:poly} we have
%
%
\begin{equation}
\corr(W_1(s_n),W_2(t_m))\to0\qquad \mbox{as $m,n\to\infty$.}
\end{equation}
\end{Proposition}

Our next theorem is the analogue of Theorem \ref{th:poly} in the
case of an exponential decay in the dependence condition.
\begin{Theorem}\label{th:exp}
Let $p>2$ and let $\{Y_k,k\in\mathbb{Z}\}$ be a centered stationary
sequence, weakly $\mathcal M$-dependent in $L^p$ with rate function
%
%
\begin{equation}\label{eq:expm}
\delta(m)\ll\exp(-\varrho m),\qquad \varrho>0.
\end{equation}
Then the series (\ref{eq:sig}) is absolutely convergent, and $\{
Y_k,k\in
\mathbb{Z}\}$ can be
redefined on a new probability space together with two standard
Wiener processes $\{W_1(t), t\geq0\}$ and $\{W_2(t), t\geq0\}$
such that
%
%
\begin{equation}\label{eq:strappr_exp}
\sum_{k=1}^n Y_k=W_1(s_n^2)+W_2(t_n^2)
+O(n^{1/p}\log^2 n) \qquad\mbox{a.s.,}
\end{equation}
where $\{s_n\}$ and $\{t_n\}$ are nondecreasing numerical sequences
such that
$s_n^2 \sim\sigma^2 n$, $t_n^2 \sim\sigma^2 n/\log n$ and (\ref
{stlimsup}) holds.
\end{Theorem}

Like in Theorem \ref{th:poly}, $s_n^2 \sim\sigma^2 n$ can be
sharpened to $s_n^2 = \sigma^2 n+O(n/\log n)$; see the remarks after
Theorem \ref{th:poly}.

Using the law of the iterated logarithm for $W_2$, relation (\ref
{eq:strappr}) implies
%
%
\begin{equation}\label{eq:strappr2}
\sum_{k=1}^n Y_k=W_1(s_n^2)+O(n^{1/2-\lambda})
\qquad\mbox{a.s.}
\end{equation}
for some $\lambda>0$, which is the standard form of strong
invariance principles. However, since $\gamma$ in (\ref{eq:sn}) is
typically near to 1,
the $\lambda$ in (\ref{eq:strappr2}) can be very small, and thus the
effect of the very strong error term $O(n^{(1+\eta)/p})$ in (\ref
{eq:strappr}) is lost.

The proof of the strong approximation theorems in Wu \cite{wuap}
depends on martingale approximation, while Liu and Lin \cite{liulin}
use approximation of the partial sums of $\{Y_k\}$ by partial sums of
$m$-dependent r.v.'s. Our approach differs from both, using a direct
approximation of separated block sums of~$\{Y_k\}$ by independent
r.v.'s, an idea used earlier in \cite{abh,behosch,berkeshorvath01,hoermann2008}.
In this approach, the second
Wiener process $W_2$ is provided by the sum of short block sums.
The question if one can get a remainder term near $o(n^{1/p})$ in the
simple (one-term) Wiener approximation for any $p>2$ remains open.




\section{Examples of weakly {$\mathcal M$}-dependent processes}\label{se:ex}


The classical approach to weak dependence, developed in the seminal
papers of Rosenblatt \cite{rosenblatt1956} and Ibragimov \cite{ibr},
uses the strong mixing property and its variants like~$\beta$,
$\varrho$, $\phi$ and $\psi$ mixing, combined with a blocking
technique
to connect the partial sum behavior of $\{Y_k\}$ with that of
independent random variables. This method yields very sharp
results (for a complete account of the classical theory see Bradley
\cite{bradley}), but
verifying mixing conditions of the above type is not easy and even when
they apply (e.g., for Markov
processes), they typically require strong smoothness conditions on
the process. For example, for the AR(1) process
\[
Y_k=\tfrac{1}{2} Y_{k-1}+\varepsilon_k
\]
with Bernoulli innovations, strong mixing fails to hold (cf.
Andrews \cite{and}). Recognizing this fact, an important line of
research in probability theory in past years has been to find
weak\vadjust{\goodbreak}
dependence conditions which are strong enough to imply satisfactory
asymptotic results, but which are sufficiently general to be
satisfied in typical applications. Several conditions of this kind
have been found, in particular by the French school (see
\cite{ded,depr1,depr2,doukhanlouhichi,rio1,rio2}).
A different type of mixing conditions, the so-called physical and
predictive dependence measures, have been introduced by
Wu \cite{wu05} for stationary processes $\{Y_k\}$ admitting the
representation
(\ref{eq:Y_k_represent}) where $\{\varepsilon_k, k\in\mathbb{Z}\}$ is
an i.i.d. sequence,
and $f\dvtx\mathbb{R}^{\mathbb N}\to\mathbb{R}$ is a Borel-measurable
function. These conditions are particularly easy to handle, since
they are defined in terms of the algorithms which generate the
process $\{Y_k\}$. Weak ${\ld}$-dependence, although formally not
requiring a representation of the form~(\ref{eq:Y_k_represent}),
is closely related to Wu's mixing conditions and works best for
processes~$\{Y_k\}$ having a representation (\ref{eq:Y_k_represent})
or its two-sided version (\ref{eq:Y_k_represent2}).
The examples below will clear up the exact connection of our weak ${\ld
}$-dependence condition with the mixing conditions in Wu \cite{wu05,wuap}
and Liu and Lin \cite{liulin}.

\subsection{$m$-dependent processes}\label{ss:mdep}
Definition \ref{d} implies that $\{Y_k, k\in\mathbb{Z}\}$ can be
approximated, for every $m\ge1$, by an $m$-dependent
process with termwise~$L^p$ error $\delta(m)$. If $\{Y_k, k\in
\mathbb
{Z}\}$ itself is $m$-dependent for some fixed $m=m_0$ and $K:={\sup
_{k\in
\mathbb{Z}}}\|Y_k\|_p<\infty$, then Definition \ref{d} is satisfied with
\[
\delta(j)=\cases{
K, &\quad if $j< m_0$,\cr
0, &\quad if $j\geq m_0$,}
\]
and $Y_k^{(n)}=0$ if $n<m_0$ and $Y_k^{(n)}=Y_k$ if $n\geq m_0$.
In other words, $m$-dependent sequences with uniformly bounded $L^p$
norms are weakly $\mathcal M$-dependent
with the above parameters.
%
%
It is worth mentioning that $m$-dependent processes in general do not have
the representation (\ref{eq:Y_k_represent2}) (see, e.g., \cite{burtonetal,valk}).

\subsection{NED processes}\label{ss:ned}
Under (\ref{eq:Y_k_represent2}) our condition can be directly compared
to NED. We recall:
\begin{Definition}[(NED)]\label{ned}
A sequence $\{Y_k, k \in{\mathbb Z}\}$ having representation (\ref
{eq:Y_k_represent2}) is called \textit{NED over $\{\varepsilon_k\}$ under
$L^p$-norm with rate function $\delta(\cdot)$} if for any $k \in
{\mathbb Z}$, $m\ge1$,
\[
\bigl\|Y_k-E[Y_k|\mathcal{F}_{k-m}^{k+m}]\bigr\|_p\leq
\delta(m),
\]
where $\mathcal{F}_{k-m}^{k+m}$ is the $\sigma$-algebra generated by
$\varepsilon_{k-m},\ldots,\varepsilon_{k+m}$.
\end{Definition}

Clearly, if $\{\varepsilon_k\}$
is an independent sequence, then $Y_k^{(m)}=E[Y_k|\mathcal
{F}_{k-m}^{k+m}]$
satisfies (B) of Definition \ref{d}. Hence if
$\{Y_k\}$ is NED over $\{\varepsilon_k\}$ in $L^p$-norm with rate
function $\delta(\cdot)$ where $\{\varepsilon_k\}$ is an independent
sequence, then $\{Y_k\}$ is weakly $\mathcal M$-dependent with the same
$p$, $\delta(\cdot)$.

As our examples below will show, for weakly ${\ld}$-dependent sequences
the construction for $Y_k^{(m)}$ is not restricted to
$E[Y_k|\mathcal{F}_{k-m}^{k+m}]$,
but is often more conveniently established by truncation or coupling methods.

\subsection{Linear processes}
\label{ss:linproc}
Let $Y_k=\sum_{j=-\infty}^{\infty}
a_j \varepsilon_{k-j}$ with the i.i.d. innovations
$\{\varepsilon_j$, $j\in\mathbb{Z}\}$. If $a_j=0$ for $j<0$, then the
sequence $\{Y_k,k\in\mathbb{Z}\}$ is causal. Liu and Lin \cite
{liulin} and
Wang, Lin and Gulati \cite{wlg} studied strong approximations
of the partial sums with Gaussian processes (in the short- and
long-memory cases).

We define $\Ykm$ as
\[
\Ykm=\sum_{j=-\lfloor m/2\rfloor}^{\lfloor m/2\rfloor} a_j
\varepsilon_{k-j}.
\]
This directly ensures that condition (B) holds. To verify
condition (\ref{eq:ma}) we
will assume that $E|\varepsilon_0|^p<\infty$ for some $p>2$ as well
as $|a_j|\ll|j|^{-(A+1)}$ $(j\to\infty)$. Then we get, using the
Minkowski inequality,
\begin{eqnarray*}
\bigl\|Y_k-\Ykm\bigr\|_p&=&\biggl\|\sum_{|j|>m/2} a_j
\varepsilon_{k-j}\biggr\|_p\\
&\le&\sum_{|j|>m/2} \|a_j \varepsilon_{k-j}\|_p\\
&=&(E|\varepsilon_0|^p)^{1/p}
\sum_{|j|>m/2} |a_j| \ll m^{-A}.
\end{eqnarray*}
Thus if $A$ is large enough, Theorem \ref{th:poly} applies.
Obviously if $|a_j|\ll\rho^{|j|}$ with some $0<\rho<1$, then
(\ref{eq:expm}) holds, and Theorem \ref{th:exp} applies.

\subsection{Nonlinear time series}
\label{ss:irf}
Let the time series $\{Y_k,k\in\mathbb{Z}\}$ be defined by the
stochastic recurrence equation
%
%
\begin{equation}\label{eq:nonlinear_ts}
Y_k=G(Y_{k-1},\varepsilon_k),
\end{equation}
where $G$ is a measurable function, and
$\{\varepsilon_k,k\in\mathbb{Z}\}$ is an i.i.d. sequence. For
example, ARCH$(1)$ processes (see, e.g., Engle \cite{engle1982})
which play an important role
in the econometrics literature, are included in this setting.
Sufficient conditions for the existence of a stationary solution of
(\ref{eq:nonlinear_ts}) can be found in Diaconis and
Freedman \cite{diaconisfreedman}.
Note that iterating
(\ref{eq:nonlinear_ts}) yields
$Y_k=f(\ldots,\varepsilon_{k-1},\varepsilon_k)$ for some measurable
function $f$. This suggests defining the
approximating random variables $Y_k^{(m)}$ as
$\Ykm=f(\ldots,0,0,\varepsilon_{k-m},\ldots,\varepsilon_k)$.
Note, however, that this definition does not guarantee
the convergence and thus the
existence of $\Ykm$. The coupling used by Wu
\cite{wu05}, avoids this problem by defining
\[
\Ykm=f\bigl(\ldots,\varepsilon_{k-m-2}^{(k)},\varepsilon_{k-m-1}^{(k)},
\varepsilon_{k-m},\ldots,\varepsilon_k\bigr),
\]
where $\{\varepsilon_{k}^{(l)},k\in\mathbb{Z}\}$, $l=1,2,\ldots,$
are i.i.d. sequences with the same distribution as
$\{\varepsilon_k,k\in\mathbb{Z}\}$ which are independent
of each other and of the $\{\varepsilon_k,\allowbreak k\in\mathbb{Z}\}$.
These random variables satisfy condition (B).\vadjust{\goodbreak} Results
from Wu and Shao \cite{shaowu} show that under some simple technical
assumption on $G$,
\[
\bigl\|Y_k-\Ykm\bigr\|_p\ll\exp(-\rho m)
\]
holds with some $p>0$ and $\rho>0$. Thus for $p>2$,
Theorem \ref{th:exp} applies.

\subsection{Augmented GARCH sequences}
\label{ss:auggarch}
Augmented GARCH sequences were introduced by Duan \cite{duan} and turned
out to be very useful in applications in macroeconomics and finance.
The model is quite general and many popular processes are included
in its framework. Among others the well-known GARCH \cite{bollerslev},
AGARCH \cite{dge} and EGARCH model \cite{nelson} are covered.
We consider the special case of augmented GARCH($1,1$) sequences,
that is, sequences $\{Y_k,k\in\mathbb{Z}\}$ defined by
%
%
\begin{equation}\label{eq:aug_GARCH}
Y_k=\sigma_k \varepsilon_k,
\end{equation}
where the conditional variance $\sigma_k^2$ is given by
%
%
\begin{equation}\label{eq:aug_GARCH_Lambda}
\Lambda(\sigma_k^2)=c(\varepsilon_{k-1})\Lambda(\sigma_{k-1}^2)
+g(\varepsilon_{k-1}).
\end{equation}
Here $\{\varepsilon_k,k\in\mathbb{Z}\}$ is a sequence of i.i.d.
errors, and $\Lambda(x)$, $c(x)$ and $g(x)$ are real-valued
measurable functions. To solve (\ref{eq:aug_GARCH_Lambda}) for
$\sigma_k^2$ one usually assumes that $\Lambda^{-1}(x)$ exists.
Necessary and sufficient conditions for the existence of a~strictly
stationary solution of (\ref{eq:aug_GARCH}) and
(\ref{eq:aug_GARCH_Lambda}) were given by Duan \cite{duan} and
Aue, Berkes and Horv{\'a}th \cite{abh}. 
Under some technical conditions stated in H\"{o}rmann \cite{hoermann2008}
(Lem\-mas~1,~2 and Remark 2) one can show that
augmented GARCH sequences are
weakly ${\ld}$-dependent in $L^p$-norm with exponential
rate.\looseness=-1

Note that the above models have short memory; long memory models (see,
e.g., \cite{girsurgrob}) have completely different properties.

\subsection{Linear processes with dependent innovations}
\label{ss:linprocdep}
Linear processes
$
Z_k=\sum_{j=-\infty}^{\infty} a_j Y_{k-j}
$
with dependent innovations $\{Y_k\}$
have obtained considerable
interest in the financial literature. A common
example are autoregressive (AR) processes with augmented
GARCH innovations (see, e.g., \cite{lientse}).

Assume that $\{Y_k,k\in\mathbb{Z}\}$ is weakly $\mathcal M$-dependent
in $L^p$ with rate function~$\delta(\cdot)$.
In combination with the results of Section \ref{ss:linproc} one can
easily obtain conditions on $\delta$ assuring that the linear process
$\{Z_k,k\in\mathbb{Z}\}$ defined above
is also weakly $\mathcal M$-dependent in $L^p$-norm with a rate
function $\delta^*$ depending on~$(a_j)$ and $\delta$.

Strong approximation results for linear processes with
dependent errors were also obtained by Wu and Min \cite{wumin}.

\subsection{Ergodic sums}
\label{ss:f2k}
Let $f$ be a real measurable function with period 1 such that $\int
_{0}^1f(\omega) \,d\omega=0$
and $\int_{0}^1|f(\omega)|^p \,d\omega<\infty$ for some $p>2$. Set
\[
S_n(\omega)=\sum_{k=1}^n f(2^k \omega),\qquad \omega\in[0,1),\vadjust{\goodbreak}
\]
and\vspace*{1pt} $B_n^2=\int_{0}^1S_n^2(\omega) \,d\omega$.
Then $S_n$ defines a partial sum process on the probability space
$([0,1),\mathcal{B}_{[0,1)},\lambda_{[0,1)})$, where $\mathcal
{B}_{[0,1)}$ and
$\lambda_{[0,1)}$ are the Borel $\sigma$-al\-gebra and Lebesgue measure
on $[0,1)$.
The strong law of large numbers for~$f(2^k\omega)$ is a consequence of
the ergodic theorem,
for central and functional central limit theorems see Kac \cite{kac},
Ibragimov \cite{ibragimov67} and~Billingsley~\cite{billingsley68}.

Let $Y_k(\omega)=f(2^k\omega)$, and define the random variable
$\varepsilon_k(\omega)$ to be equal to the $k$th digit in the binary
expansion of $\omega$. Ambiguity can be avoided by the convention to
take terminating expansions whenever possible. Then~%
$\{\varepsilon_k\}$ is an i.i.d. sequence, and we
have $\varepsilon_k=\pm1$, each with probability $1/2$. This gives the
representation
\[
Y_k=f\Biggl(\sum_{j=1}^{\infty}\varepsilon_{k+j}2^{-j}\Biggr)=
g(\varepsilon_{k+1},\varepsilon_{k+2},\ldots).
\]
We can now make use of the coupling method described in Section
\ref{ss:irf} and the approximations
\[
\Ykm=g\bigl(\varepsilon_{k+1},\varepsilon_{k+2},\ldots,\varepsilon_{k+m},
\varepsilon_{k+m+1}^{(k)},\varepsilon_{k+m+2}^{(k)},\ldots\bigr).
\]
Changing for some $\omega\in[0,1)$ the digits
$\varepsilon_k(\omega)$ for $k> m$ will give an $\omega'$ with
$|\omega-\omega'|\leq2^{-m}$. If $f$ is Lipschitz continuous of
some order $\gamma$, then we have
\[
\bigl|Y_k-\Ykm\bigr|=O(2^{-\gamma m}),
\]
and thus for any $p\geq1$ $\{Y_k\}$ is weakly ${\ld}$-dependent in
$L^p$-norm with
an exponentially decaying rate function.


\section{Increments of stochastic processes}\label{se:app}

For arbitrary $\lambda>0$, relation (\ref{eq:strappr2}) has many useful
applications in probability and statistics. For example,
it implies a large class of limit theorems on CLT and LIL behavior and
for various other functionals of weakly dependent sequences.
However, many refined limit theorems for partial sums require a
remainder term better than~$O(n^{1/4})$, and no existing method for
dependent sequences
provides such a~remainder term. The purpose of the next two sections is
to show how to deal with such limit theorems via our approximation
results in Section \ref{se:main}.


Let $\{Y_k,k\in\mathbb{Z}\}$ be a~stationary random sequence, and let
$0<a_n\leq n$ be a~nondecreasing sequence of
real numbers. In this section, we investigate the order of magnitude of
\[
\max_{1\leq k \leq n-a_n}\max_{1\leq\ell\leq
a_n}\Biggl|\sum_{j=k+1}^{k+\ell}Y_j\Biggr|.
\]
Such results have been obtained by Cs\"{o}rg\H{o} and R\'{e}v\'{e}sz
\cite{cr} for i.i.d. sequences and the Wiener process. In particular,
they obtained the following result (\cite{cr}, Theorem~1.2.1).
%
\begin{Theorem}\label{csr_incr}
Let
$\{a_T, T\ge0\}$ be a positive nondecreasing function satisfying:
\begin{longlist}[(a)]
\item[(a)] $0<a_T\leq T$;
\item[(b)] $T/a_T$ is nondecreasing. 
\end{longlist}
Set
%
%
\begin{equation} \label{beta}
\beta_T=\biggl(2a_T\biggl[\log\frac{T}{a_T}+\log\log
T\biggr]\biggr)^{-1/2}.
\end{equation}
%
Then
\[
\mathop{\overline{\lim}}_{T\to\infty}\max_{0\leq t\leq T-a_T}\max_{0\leq s
\leq a_T}\beta_T|W(t+s)-W(t)|= 1.
\]
%
\end{Theorem}

Using strong invariance, a similar result can be obtained for partial
sums of i.i.d. random variables under suitable moment conditions (see
\cite{cr}, pa\-ges~115--118). For slowly growing $a_T$, this requires a
very good remainder term in the Wiener approximation of partial sums,
using the full power of the Koml\'{o}s--Major--Tusn\'{a}dy theorems. As
an application of our main theorems in Section~\ref{se:main}, we now
extend Theorem \ref{csr_incr} for dependent stationary processes. To
simplify the formulation and to clarify the connection between the
remainder term in our approximation theorems in Section \ref{se:main}
and the increment problem, we introduce the following
assumption.
\begin{Assumption}\label{ass:strappr} Let $\{Y_k\}$ be a
random sequence which can be redefined on a new probability space
together with two standard Wiener processes
$\{W_1(t), t\geq0\}$ and $\{W_2(t), t\geq0\}$ such that
%
%
\begin{equation}\label{eq:assstrappr}
\sum_{k=1}^n
Y_k=W_1(s_n^2)+W_2(t_n^2)+O(E_n)
\qquad\mbox{a.s.},
\end{equation}
where $\{E_n\}$ is some given sequence and $\{s_n^2\}$ and $\{t_n^2\}$
are nondecreasing sequences satisfying
%
%
\begin{equation}\label{stcond}\qquad
s_n^2\sim\sigma^2 n,\qquad t_n^2=o(n),\qquad \mathop{\overline{\lim}}_{k\to
\infty} (s_{k+1}^2-s_k^2)=\mathop{\overline{\lim}}_{k\to\infty}
(t_{k+1}^2-t_k^2)= \sigma^2.
\end{equation}
\end{Assumption}

We 
will prove the following result.
\begin{Theorem}\label{th:fluct}
Let $\{Y_k\}$ be a sequence of random variables satisfying
Assumption \ref{ass:strappr} and put $S_n=\sum_{k=1}^n Y_k$. Let $a_T$
be a positive nondecreasing function such that:
\begin{longlist}[(a)]
\item[(a)] $0<a_T\leq T$;
\item[(b)] $T/a_T$ is nondecreasing;
\item[(c)]
$a_T$ is regularly varying at $\infty$ with index $\varrho\in(0,1]$.
\end{longlist}
%
Let $\beta_T$ be defined by (\ref{beta}). Then under the condition
%
%
\begin{equation}\label{eq:gh}
\beta_T E_T=o(1)
\end{equation}
we have
%
%
\begin{equation}\label{concl}
\mathop{\overline{\lim}}_{n\to\infty}\max_{1\leq k\leq n-a_n}\max_{1\leq
\ell
\leq a_n}\beta_n|S_{k+\ell}-S_k|= \sigma^2.
\end{equation}
\end{Theorem}

Given a function $a_T$ and a weakly $\mathcal M$-dependent sequence $\{
Y_k\}$
with parameters $p, \delta(\cdot)$, we can compute, using Theorem
\ref{th:fluct}, a rate of decrease for $\delta(\cdot)$ and a value for
$p>2$ such that the fluctuation result (\ref{concl}) holds. For
example, if $a_T=T^\alpha$, $0<\alpha<1$, then (\ref{concl}) holds if
$p>4/\alpha$ and $\delta(m)\ll m^{-p/2}$.\looseness=1

We note that for i.i.d. observations
only assumptions (a) and (b) are required.
It remains open whether a more general version of our Theorem \ref
{th:fluct} which does not require assumption (c) can be proved.

Recently Zholud \cite{zholud1} obtained a distributional version of
Theorem \ref{csr_incr} by showing that the functional
\[
\max_{0\leq t\leq T-a_T}\max_{0\leq s \leq a_T}\bigl(W(t+s)-W(t)\bigr)
\]
converges weakly, suitably centered and normalized, to the extremal
distribution with distribution function $e^{-e^{-x}}$. Using this fact
and our a.s. invariance principles, a distributional version of
Theorem \ref{th:fluct} can be obtained easily. Since the argument is
similar to that for (\ref{concl}), we omit the details.\looseness=1




Let $|A|$ be the cardinality of a set $A$.
For the proof of Theorem \ref{th:fluct} we need the following simple lemma.

\begin{Lemma}\label{le:dens}
Assume that $\{d_k,k\ge1\}$ is a nonincreasing sequence of
positive numbers such that $\sum_{k=1}^{\infty} d_k=\infty$. Let
$A\subset
\mathbb{N}$ have positive density, that is,
\[
\liminf_{n\to\infty}|A\cap\{1,\ldots,n\}|/n>0.
\]
Then $\sum_{k=1}^{\infty}
d_kI\{k\in A\}=\infty$.
\end{Lemma}
\begin{pf}
First note that by our assumption
we have $\sum_{k=1}^nI\{k\in A\}\geq\mu n$ for some $\mu>0$ as long as
$n\geq n_0$.
Using Abel summation we can write
\[
\sum_{k=1}^nd_k=nd_n+\sum_{k=1}^{n-1}k(d_k-d_{k+1}).
\]
Hence, by our assumptions
\[
nd_n+\sum_{k=n_{0}}^{n-1}k(d_k-d_{k+1})\to\infty\qquad (n\to\infty).
\]
From $d_k-d_{k+1}\geq0$ it follows (again using the Abel summation) that
for $n\geq n_0$
\begin{eqnarray*}
\hspace*{50pt}&&\sum_{k=1}^nd_kI\{k\in A\}\\[-2pt]
&&\qquad= d_n\sum_{k=1}^nI\{k\in
A\}+\sum_{k=1}^{n-1}(d_k-d_{k+1})\sum_{j=1}^kI\{j\in
A\}\\[-2pt]
&&\qquad\geq\mu\Biggl(
nd_n+\sum_{k=n_0}^{n-1}k(d_k-d_{k+1})\Biggr)\to\infty\qquad (n\to
\infty).\hspace*{38pt}\qed
\end{eqnarray*}
\noqed\end{pf}
\begin{pf*}{Proof of Theorem \ref{th:fluct}}
For the sake of simplicity we carry out the proof for $\sigma=1$.
From (\ref{eq:assstrappr}) and the triangular inequality we infer
that
\begin{eqnarray*}
&&\mathop{\overline{\lim}}_{n\to\infty}\max_{1\leq k\leq n-a_n}\max_{1\leq
\ell
\leq
a_n}\beta_n|S_{k+\ell}-S_k|\\[-2pt]
&&\qquad \leq\mathop{\overline{\lim}}_{n\to\infty}\max_{1\leq k\leq
n-a_n}\max
_{1\leq\ell\leq
a_n}\beta_n|W_1(s_{k+\ell}^2)-W_1(s_k^2)
|\\[-2pt]
&&\qquad\quad{} +\mathop{\overline{\lim}}_{n\to\infty}\max_{1\leq k\leq n-a_n}\max
_{1\leq
\ell\leq
a_n}\beta_n|W_2(t_{k+\ell}^2)-W_2(t_k^2)
|\\[-2pt]
&&\qquad\quad{} +\mathop{\overline{\lim}}_{n\to\infty}\beta_n O(E_n)\\[-2pt]
&&\qquad = A_1+A_2+A_3.
\end{eqnarray*}
By (\ref{eq:gh}) $A_3=0$. Since $a_n\to\infty$ [this is implicit in
(c)], we conclude from (\ref{stcond}) that for any
$\varepsilon>0$ some $n_0$ exists, such that for all $n\geq
n_0$
\[
\sup_{k\geq1} \{s_{k+a_n}^2-s_{k}^2\}\leq
(1+\varepsilon)a_n \quad\mbox{and}\quad s_n^2\leq
(1+\varepsilon)n.
\]
Set $T=(1+\varepsilon)n$, and define
$
a_{T,\varepsilon}=(1+\varepsilon)a_{T/(1+\varepsilon)}.
$
Then $a_{T,\varepsilon}$ satisfies (a) and (b)
%
%
and for $n\geq n_0$ we have
\begin{eqnarray*}
&&\max_{1\leq k\leq n-a_n}\max_{1\leq\ell\leq
a_n}\beta_n|W_1(s_{k+\ell}^2)-W_1(s_k^2)
|\\[-2pt]
&&\qquad \leq\sup_{0\leq t\leq s^2_{n-a_n}}\sup_{0\leq s\leq
(1+\varepsilon)a_n}
\beta_n|W_1(t+s)-W_1(t)|\\[-2pt]
&&\qquad \leq\sup_{0\leq t\leq T-a_{T,\varepsilon}}\sup_{0\leq s\leq
a_{T,\varepsilon}}
\beta_{T/(1+\varepsilon)}|W_1(t+s)-W_1(t
)|.
\end{eqnarray*}
Let
\[
\beta_{T,\varepsilon}=\biggl(2a_{T,\varepsilon}\biggl[\log\frac
{T}{a_{T,\varepsilon}}+
\log\log T\biggr]\biggr)^{-1/2}.
\]
By application of Theorem 1.2.1 in Cs\"{o}rg\H{o} and R\'{e}v\'{e}sz
\cite{cr} [which requires~(a) and (b)] we get
\[
\mathop{\overline{\lim}}_{T \to\infty}\sup_{0\leq t\leq T-a_{T,\varepsilon
}}\sup
_{0\leq s\leq
a_{T,\varepsilon}}
\beta_{T,\varepsilon}|W_1(t+s)-W_1(t)
|=1
\qquad\mbox{a.s.}\vadjust{\goodbreak}
\]
%
Since
$
\overline{\lim}_{T\to\infty}\beta_{T/(1+\varepsilon)}/\beta
_{T,\varepsilon}= (1+\varepsilon)^{1/2}
$,
and $\varepsilon$ can be chosen arbitrarily small, we have shown that
$A_1\leq1$ a.s.

It is not surprising that due to (\ref{stcond})
similar arguments will lead to $A_2=0$ a.s.

The proof will be completed if we show that $A_1\geq1$. Let $\{n_k\}$ be
a nondecreasing sequence of integers with $n_k\to\infty$.
By
(\ref{eq:assstrappr}), the triangular inequality\vadjust{\goodbreak} and $A_2=A_3=0$ we obtain
\begin{eqnarray*}
&&\mathop{\overline{\lim}}_{n\to\infty}\max_{1\leq k\leq n-a_n}\max_{1\leq
\ell
\leq
a_n}\beta_n|S_{k+\ell}-S_k|\\
&&\qquad\geq \mathop{\overline{\lim}}_{n\to\infty}\max_{1\leq k\leq n-a_n}\max
_{1\leq
\ell\leq
a_n}\beta_n|W_1(s_{k+\ell}^2)-W_1(s_k^2)
|\\
&&\qquad\geq \mathop{\overline{\lim}}_{k\to\infty}\beta_{n_k}|W_1
(s_{n_k}^2)-W_1
(s_{n_k-a_{n_k}}^2)|.
\end{eqnarray*}
We now proceed similarly as in Cs\"{o}rg\H{o} and
R\'{e}v\'{e}sz \cite{cr} for the proof of Step 2 of their Theorem 1.2.1.
We will distinguish between the cases $\lim
a_T/T=\rho$ with $\rho<1$ and $\rho=1$. Since both times we can
use the same conceptual idea, we shall treat here only
$\rho<1$.

Set $n_1=1$. Given $n_k$, define $n_{k+1}$ such that
$n_{k+1}-a_{n_{k+1}}=n_k$. This equation will, in general,
have no integer solutions, but for the sake of simplicity we assume that
$(n_k)$ and $(a_{n_k})$ are $\mathbb{Z}$-valued. Since
$(s_n^2)$ is nondecreasing, we conclude that the
increments
$
\Delta(k)=W_1(s_{n_k}^2)-W_1(s_{n_{k-1}}^2)
$
are independent. By the second Borel--Cantelli lemma it suffices
to show now that
%
%
\begin{equation}\label{eq:bc}
\sum_{k=1}^{\infty}P\bigl(\beta_{n_k}|\Delta(k)|\geq
1-\varepsilon\bigr)=\infty\qquad\mbox{for all $\varepsilon>0$.}
\end{equation}
%
For all large enough $k\in\mathbb{N}$ for which
$
s_{n_k}^2-s_{n_{k-1}}^2\geq(1-\varepsilon/2)a_{n_k},
$
the estimates in \cite{cr} give
\[
P\bigl(\beta_{n_k}|\Delta(k)|\geq1-\varepsilon\bigr)\geq
\biggl(\frac{a_{n_k}}{n_k\log n_k}\biggr)^{1-\varepsilon}.
\]
It is also shown in \cite{cr}
that $\sum_{k=1}^{\infty}(\frac{a_{n_k}}{n_k\log n_k}
)^{1-\varepsilon}=\infty$.
Thus, in view of Lemma \ref{le:dens} it remains to show that $A=\{
k\geq1|
s_{n_k}^2-s_{n_{k-1}}^2\geq(1-\varepsilon/2)a_{n_k}\}$ has
a positive density. By (\ref{stcond}) we have
\begin{eqnarray*}
(s_{n_k}^2-s_{n_1}^2)/n_k&=&\sum_{j=2}^k(s_{n_j}^2-s_{n_{j-1}}^2)/n_k\\
&\leq& C_0\mathop{\sum_{2\leq j\leq k}}_{j\in A}(n_j-n_{j-1})/n_k+
\mathop{\sum_{2\leq j\leq k}}_{j\in A^c}(1-\varepsilon/2)(n_j-n_{j-1})/n_k\\
&\leq& C_0\mathop{\sum_{2\leq j\leq k}}_{j\in
A}(n_j-n_{j-1})/n_k+(1-n_1/n_k)(1-\varepsilon/2)
\end{eqnarray*}
for some $C_0>0$ which is independent of $k$.
Now if $A$ had density zero, the limsup of the right-hand side
of the last relation would be $1-\varepsilon/2$.
This can be easily proved, using that $(n_j-n_{j-1})$ is regularly varying
by assumption~(c).
The liminf of the left-hand side above is 1.
Thus $A$ must have positive density and the proof is complete.
\end{pf*}

\section{Change-point tests with an epidemic alternative}\label{se:epi}

In this section we apply our invariance principles to a
change-point problem.
Let $\{Y_k,k\in\mathbb{Z}\}$ be a zero mean process.
Further let $X_k=Y_k+\mu_k$, where
$\mu_k$, $k\in\mathbb{Z}$, are unknown constants. We want to
test the hypothesis
{\renewcommand{\theequation}{$H_0$}
\begin{equation}\label{equH0}
\mu_1=\mu_2=\cdots=\mu_n=\mu
\end{equation}}

\vspace*{-\baselineskip}

\noindent against the ``epidemic alternative''
{\renewcommand{\theequation}{$H_A$}
\begin{equation}\label{equHa}
\begin{tabular}{p{315pt}}
\mbox{There exist $1\le m_1< m_2\le n$ such that
$\mu_k=\mu$ for}
\mbox{$k\in\{1,\ldots,
n\}\setminus\{m_1+1, \ldots, m_2\}$ and
$\mu_k=\mu+\Delta$} \mbox{if $k\in\{m_1+1,\ldots,m_2\}$.}
\end{tabular}\hspace*{-30pt}
\end{equation}}

\vspace*{-12pt}

\noindent It should be noted that the variables $m_1$, $m_2$ and $\Delta$ may
depend on the sample size~$n$. As it is common in the change-point
literature, this dependence is suppressed in the notation.

Without loss of generality we assume that $\sigma=1$. To detect a possible
epidemic change it is natural to compare the increments
of the process to a proportion of
the total sum. More specifically, assume for the moment
that $X_k$ are independent and that we know when the
epidemic starts and ends. Set $S_k=X_1+\cdots+X_k$. Then
by the law of large numbers
$I(m_1,m_2)=|S_{m_2}-S_{m_1} - (m_2-m_1) S_n/n|\gg m_2-m_1$.
If no change occurs, however, by the central limit theorem
$I(m_1,m_2)=O_P(\sqrt{m_2-m_1})$.
In general we do not know
$m_1$ and $m_2$. Thus, a natural test statistic is
\[
\max_{1\le i < j \le n} |S_j-S_i-(j-i) S_n/n|.
\]
Clearly
we are required to normalize the above
test statistic appropriately. Following Ra\u{c}kauskas and
Suquet \cite{rasu} we define
\[
UI(n,\alpha)=n^{-1/2} \max_{1\le i < j \le n}
\frac{|S_j-S_i-(j-i) S_n/n|}{[((j-i)/n)
(1-(j-i)/n)]^{\alpha}}
\]
with $0<\alpha<1/2$.
As we will see below, the parameter
$\alpha$ plays an important role. The closer $\alpha$ is
to $1/2$, the ``shorter'' epidemics can be detected with this
test. The price, however, is that in order
to obtain the limiting law under (\ref{equH0}) with ``large'' $\alpha$
(close to $1/2$)
requires a.s. invariance principles with error $n^{\epsilon}$,
$\epsilon$
close to zero.
Choosing $\alpha\geq1/2$ would result in a degenerate
limiting distribution under (\ref{equH0}).
%
%
\begin{Proposition}[{[Asymptotics under (\ref{equH0})]}]\label{pr:ecH0}
If the\vspace*{1pt} stationary sequence $\{Y_k, k\in\mathbb{Z}\}$ satisfies
Assumption \ref{ass:strappr} with $E_n=o(n^{1/2-\alpha})$ and
(\ref{equH0}) holds, then
\[
\sigma^{-1} UI(n,\alpha) \convd\sup_{0<s<t<1}
\frac{|B(t)-B(s)|}{[(t-s)(1-(t-s))]^{\alpha}},
\]
where
$\{B(t),t\in[0,1]\}$ is a Brownian bridge.
\end{Proposition}
\begin{pf}
Using (\ref{eq:assstrappr}) and assuming for simplicity that
$\sigma=1$, we obtain
\begin{eqnarray*}
UI(n,\alpha) &\le& n^{-1/2}\max_{1\le i<j\le n}
\frac{|W_1(s_j^2)-W_1(s_i^2)-(j-i) W_1(s_n^2)/n|}{[((j-i)/n)
(1-(j-i)/n)]^{\alpha}}\\[-0.8pt]
&&{}+ n^{-1/2}\max_{1\le i<j\le n}
\frac{|W_2(t_j^2)-W_2(t_i^2)-(j-i) W_2(t_n^2)/n|}{[((j-i)/n)
(1-(j-i)/n)]^{\alpha}}\\[-0.8pt]
&&{}+ O(n^{-1/2+\alpha}E_n)\\[-0.8pt]
&=&n^{-1/2}\max_{1\le i<j\le n}T_{i,j}^{(1)}+n^{-1/2}\max_{1\le i<j\le
n}T_{i,j}^{(2)}+o(1).
\end{eqnarray*}
It is easy to
see that $n^{-1/2}\max_{1\le i<j\le n}T_{i,j}^{(2)}$ tends to zero.
Since we can get a similar lower bound for $UI(n,\alpha)$, we have
\[
UI(n,\alpha)=n^{-1/2}\max_{(i,j)\in\mathcal{M}_n}T_{i,j}^{(1)}+o_P(1),
\]
where
$\mathcal{M}_n=\{(i,j)| 1\leq i<j\leq n\}$. Let us partition $\mathcal
{M}_n$ into
\begin{eqnarray*}
\mathcal{M}_{1,n}&=&\{(i,j)| 1\leq i<j\leq n; n\gamma
_n<j-i<n(1-\gamma
_n)\},\\[-0.8pt]
\mathcal{M}_{2,n}&=&\{(i,j)| 1\leq i<j\leq n; n\gamma_n\geq j-i\}
\end{eqnarray*}
and
\[
\mathcal{M}_{3,n}=\{(i,j)| 1\leq i<j\leq n; j-i\geq
n(1-\gamma_n)\},
\]
where $\gamma_n\to0$ will be defined
later. By our assumptions on the sequence $\{s_j^2\}$ there exists a
$\tau>0$ such that $s_j^2-s_i^2\leq\tau(j-i)$ for all\vspace*{2pt}
$1\leq i\le j$ and that $s_n^2\leq(2-\tau\gamma_n)n$ if $n\geq n_0$. We
have for large enough $n$
\begin{eqnarray*}
&&n^{-1/2}\max_{(i,j)\in\mathcal{M}_{2,n}}T_{i,j}^{(1)}\\[-0.8pt]
&&\qquad \leq2n^{\alpha-1/2}\max_{(i,j)\in\mathcal{M}_{2,n}}\biggl\{
\frac{|W_1(s_j^2)-W_1(s_i^2)|}{(j-i)^{\alpha}}\biggr\}+2n^{-1/2}
\gamma_n^{1-\alpha}|W_1(s_n^2)|\\[-0.8pt]
&&\qquad \leq2n^{\alpha-1/2}\max_{1\leq h\leq n\gamma_n}\sup_{0\leq
t\leq
(2-\tau\gamma_n)n}\sup_{0\leq s\leq\tau h}\biggl\{
\frac{|W_1(t+s)-W_1(t)|}{h^{\alpha}}\biggr\}+o_P(1).
\end{eqnarray*}
%
For arbitrary $\epsilon>0$ we get by Lemma 1.2.1 in Cs\"{o}rg\H{o} and
R\'{e}v\'{e}sz \cite{cr}
that there is a constant $C$ which is independent of $n$ and $\epsilon$
such that
\begin{eqnarray*}
&&P\biggl(\max_{1\leq h\leq n\gamma_n}\sup_{0\leq t\leq(2-\tau
\gamma
_n)n}\sup_{0\leq s\leq\tau h}\biggl\{
\frac{|W_1(t+s)-W_1(t)|}{h^{\alpha}}\biggr\}>\epsilon
n^{1/2-\alpha
}\biggr)\\[-0.8pt]
&&\qquad \leq\sum_{h=1}^{n\delta_n}
P\Bigl(\sup_{0\leq t\leq2n-\tau h}\sup_{0\leq s\leq\tau h}
|W_1(t+s)-W_1(t)|>\epsilon h^{1/2} (n/h)^{1/2-\alpha}\Bigr)\\[-0.8pt]
&&\qquad \leq\sum_{h=1}^{n\delta_n}\frac{Cn}{h} e^{-({\epsilon
^2}/{3}) (n/h)^{1-2\alpha}}\to0\qquad (n\to\infty).
\end{eqnarray*}
Hence $n^{-1/2}\max_{(i,j)\in\mathcal{M}_{2,n}}T_{i,j}^{(1)}=o_P(1)$.
In the same fashion one can show that
$n^{-1/2}\max_{(i,j)\in\mathcal{M}_{3,n}}$ $T_{i,j}^{(1)}=o_P(1)$.
Therefore
\[
UI(n,\alpha)=n^{-1/2}\max_{(i,j)\in\mathcal{M}_{1,n}}T_{i,j}^{(1)}+o_P(1).
\]
Some further basic estimates give
\begin{eqnarray*}
&&n^{-1/2}\max_{(i,j)\in\mathcal{M}_{1,n}}T_{i,j}^{(1)}\\
&&\qquad = n^{-1/2}\max_{(i,j)\in\mathcal{M}_{1,n}}\frac
{|W_1(j)-W_1(i)-(j-i)W_1(n)/n|}
{[((j-i)/n) (1-(j-i)/n)]^{\alpha}}\\
&&\qquad\quad{} +O\biggl(\frac{n^{-1/2}}{\gamma_n^{\alpha}}
\max_{1\leq i\leq n}|W_1(i)-W_1(s_i^2)|\biggr).
\end{eqnarray*}
Since $s_n^2\sim n$ there is a null sequence
$\{\epsilon_n\}$ such that $\max_{1\leq i\leq n}|i-s_i^2|\leq
\epsilon
_n n$.
Hence
\[
\max_{1\leq i\leq n}|W_1(i)-W_1(s_i^2)|\leq\sup_{0\leq t \leq n}\sup
_{0\leq s\leq2\epsilon_n n}
|W_1(t+s)-W_1(t)|.
\]
Setting $\gamma_n=\epsilon_n$ and applying
again Lemma 1.2.1 in \cite{cr} it can be seen that
\[
\sup_{0\leq t \leq n}\sup_{0\leq s\leq2\epsilon_n n}
|W_1(t+s)-W_1(t)|=o_P(n^{1/2}\gamma_n^{\alpha}).
\]
Consequently
%
%
\setcounter{equation}{27}
\begin{eqnarray}\label{eq:UI}
UI(n,\alpha)&=&n^{-1/2}\max_{(i,j)\in\mathcal{M}_{1,n}}\frac
{|W_1(j)-W_1(i)-(j-i)W_1(n)/n|}
{[((j-i)/n) (1-(j-i)/n)]^{\alpha}}\nonumber\\[-8pt]\\[-8pt]
&&{}+o_P(1).\nonumber
\end{eqnarray}
Since the line of argumentation is very similar to what we
have shown before, we note now without proof
that $\mathcal{M}_{1,n}$ in the right-hand side of (\ref{eq:UI})
can be replaced by $\mathcal{M}_n$.
The rest of the proof of Proposition \ref{pr:ecH0} is standard.
\end{pf}

The next proposition shows that this test is consistent. Let
$\ell=m_2-m_1$ denote the length of the epidemic.
\begin{Proposition}[{[Asymptotics under (\ref{equHa})]}]\label{pr:ecH1}
Let $\{Y_k,k\in\mathbb{Z}\}$ be a mean zero process, weakly $\mathcal
M$-dependent in $L^p$ with $p\geq2$ and $\delta(\cdot)$ satisfying
\[
\sum_{m\geq1}\delta(m)<\infty.
\]
Let $X_k=Y_k+\mu_k$, $k\in\mathbb{Z}$. Assume that (\ref{equHa}) holds and
that
%
%
\begin{equation}\label{eq:H1_consistency}
\lim_{n\to\infty}
\frac{(\ell(n-\ell))^{1-\alpha}}{n^{3/2-2\alpha}} |\Delta
|=\infty.
\end{equation}
Then $UI(n,\alpha)\convP\infty$.
\end{Proposition}
%
%
\begin{pf}
Under the alternative hypothesis (\ref{equHa}) we have $X_k=Y_k+\mu$ for
$k\in\{1,\ldots,n\}\setminus\{m_1 +1,\ldots,m_2 \}$ and
$X_k=Y_k+\mu+\Delta$ for $k\in\{m_1 +1,\ldots,m_2 \}$. To find a
lower bound for $UI(n,\alpha)$ we study the numerator of the test
statistic corresponding to the true epidemic. Thus we look at
\begin{eqnarray*}
&&S_{m_2 }-S_{m_1 }-S_n(m_2 /n-m_1 /n)\\
&&\qquad=(1-\ell/n)
(S_{m_2 }-S_{m_1 })-(\ell/n)\bigl(S_n-(S_{m_2 }-S_{m_1 })\bigr)\\
&&\qquad=\frac{\ell(n-\ell)}{n} \Delta+ (1-\ell/n) \sum_{j=m_1 +1}^{m_2 }
Y_j-(\ell/n)\Biggl(\sum_{j=1}^{m_1 } Y_j+\sum_{j=m_2 +1}^n Y_j
\Biggr)\\
&&\qquad=\frac{\ell(n-\ell)}{n} \Delta+ R_n.
\end{eqnarray*}
With the help of the moment inequality stated in Proposition \ref
{pr:mi} below we get
\begin{eqnarray*}\Var(n^{-1/2}R_n)&=&O
\bigl((1-\ell/n)^2(\ell/n)+(\ell/n)^2(1-\ell/n)\\
&&\hspace*{60.6pt}{}+2(\ell/n)^{3/2}(1-\ell
/n)^{3/2}\bigr)\\
&=&O\bigl((\ell/n)(1-\ell/n)\bigr),
\end{eqnarray*}
and thus
$n^{-1/2}R_n=O_P((\ell/n)^{1/2}(1-\ell/n)^{1/2})$. Thus we
have shown that
%
%
\begin{eqnarray}\label{eq:H1_lower}
UI(n,\alpha)&\ge& n^{1/2}
\bigl((\ell/n)(1-\ell/n)\bigr)^{1-\alpha} |\Delta|\nonumber\\
&&{}-
O_P\bigl(\bigl((\ell/n)(1-\ell/n)\bigr)^{1/2-\alpha}\bigr)
\nonumber\\[-8pt]\\[-8pt]
&=&\frac{(\ell(n-\ell))^{1-\alpha}}{n^{3/2-2\alpha}} |\Delta|\nonumber\\
&&{}-O_P\bigl(\bigl((\ell/n)(1-\ell/n)\bigr)^{1/2-\alpha}\bigr).
\nonumber
\end{eqnarray}
To conclude\vspace*{1pt} the proof we note that $\lim_{n\to\infty}
((\ell/n)(1-\ell/n))^{1/2-\alpha}=0$ if $\ell=\allowbreak o(n)$ [or
$n-\ell=o(n)$, resp.] and
$((\ell/n)(1-\ell/n))^{1/2-\alpha}\le1$ in general.
Consequently condition (\ref{eq:H1_consistency}) together with
relation (\ref{eq:H1_lower}) finishes the proof.
\end{pf}

For example, if $\Delta$ is independent of $n$, then condition (\ref
{eq:H1_consistency}) will hold for
$\ell\sim c n$, $c\in(0,1)$. In case that $n^\nu\ll\ell\ll
n-n^{\nu
}$, $\nu>0$, 
condition (\ref{eq:H1_consistency}) holds
provided that $(1-2\alpha)/(1-\alpha)<2\nu$. That is, choosing $\alpha$ close to
$1/2$ allows us to detect relatively ``short'' (``long'') epidemics.

\section{Proof of the main theorems}\label{se:proofs}

\subsection{A moment inequality}
In the proofs of our theorems we will use the following moment
inequality which may be of separate interest.
%
\begin{Proposition}\label{pr:mi}
Let $\{Y_k, k\in\mathbb{Z}\}$ be a centered stationary sequence,
weak\-ly $\mathcal M$-dependent in $L^p$ with $p\geq2$ and a rate
function $\delta(\cdot)$ satisfying
\[
D_p:=\sum_{m=0}^{\infty}\delta(m)<\infty.
\]
%
Then for any $n\in\mathbb{N}$, $b\in\mathbb{Z}$ we have
%
%
\begin{equation}\label{eq:ri}
E\Biggl|\sum_{k=b+1}^{b+n}Y_k\Biggr|^p\leq C_p n^{p/2},
\end{equation}
where $C_p$ is a constant depending on $p$ and the sequence $\{Y_k\}$.
\end{Proposition}
\begin{pf}
By stationarity, we can assume $b=0$. Let first $p=2$. We use below that
$\sup_{m\geq0}\|Y_k^{(m)}\|_p\leq\|Y_1\|_p+D_p$.
Without loss of generality we assume that $EY_k^{(m)}=0$ for all $k\in
\mathbb{Z}$ and $m\in\mathbb{N}$. Since
\begin{eqnarray*}
Y_kY_{k+j}&=&\bigl(Y_k-Y_k^{(j-1)}\bigr)Y_{k+j}+Y_k^{(j-1)}
\bigl(Y_{k+j}-Y_{k+j}^{(j-1)}\bigr)\\
&&{}+Y_k^{(j-1)}Y_{k+j}^{(j-1)},
\end{eqnarray*}
we get by assumption (B) that for $j\geq1$
%
%
\begin{eqnarray}\label{eq:sighelp}
|EY_kY_{k+j}|&\leq&\bigl|E\bigl[\bigl(Y_k-Y_k^{(j-1)}
\bigr)Y_{k+j}\bigr]\bigr|
+\bigl|E\bigl[Y_{k}^{(j-1)}\bigl(Y_{k+j}-Y_{k+j}^{(j-1)}\bigr)\bigr]
\bigr|\nonumber\\
&\leq&\|Y_{k+j}\|_2 \bigl\|Y_k-Y_k^{(j-1)}\bigr\|_2 +
\bigl\|Y_k^{(j-1)}\bigr\|_2 \bigl\|Y_{k+j}-Y_{k+j}^{(j-1)}\bigr\|_2
\nonumber\\[-8pt]\\[-8pt]
&\leq&\bigl(\|Y_{k+j}\|_2+\bigl\|Y_k^{(j-1)}\bigr\|_2
\bigr)\delta(j-1)
\nonumber\\
&\leq&(2\|Y_1\|_2+D_2)\delta(j-1).\nonumber
\end{eqnarray}
From relation (\ref{eq:sighelp}) we infer, letting $S_n=\sum_{k=1}^n Y_k$,
\begin{eqnarray*}
ES_n^2&=&\sum_{k=1}^n EY_k^2+2\sum_{1\leq k<l\leq n}EY_kY_l\\
&\leq& n \|Y_1\|_2^2+2\biggl[\sum_{1\leq k\leq
n-1}|EY_kY_{k+1}|+\cdots+\sum_{1\leq k\leq2
}|EY_kY_{k+n-2}|+E|Y_1Y_n|\biggr]\\
&\leq& n\|Y_1\|_2^2+ 2(2\|Y_1\|_2+D_2)[ (n-1)\delta
(0)+\cdots+2
\delta(n-3)+\delta(n-2)]\\
&\leq& n\bigl(\|Y_1\|_2^2 +2D_2(2\|Y_1\|_2+D_2)\bigr)=:C_2 n.
\end{eqnarray*}
This shows (\ref{eq:ri}) for $p=2$.\vadjust{\goodbreak}

Once (\ref{eq:ri}) is established for $p$, it holds for all $0<q\leq p$.
Indeed, by Lyapunov's inequality, relation (\ref{eq:ri}) implies
%
%
\begin{equation}\label{eq:riq}
E\Biggl|\sum_{k=b+1}^{b+n}Y_k\Biggr|^q\leq C_p^{q/p}n^{q/2}
\end{equation}
for any $0<q\leq p$. In particular, (\ref{eq:ri}) holds with
$p=1$.

Next we prove (\ref{eq:ri}) for all integers $p>2$.
Clearly, if $C_p\geq\|Y_1\|_p^p$, then the inequality
%
%
\begin{equation}\label{eq:ri0}
E|S_n|^p \leq C_p n^{p/2}
\end{equation}
holds for $n=1$. Using a double induction argument, we show now that
for some constant $C_p$, relation (\ref{eq:ri0}) holds for all $n\in
{\mathbb N}$. More precisely, we show that if (\ref{eq:ri0}) holds for
$p-1$ and all $n\in{\mathbb N}$ and also for $p$ and $n\leq n_0$, then
it will also hold for $p$ and $n\leq2n_0$.

For $k\leq n$ put $S_{k}^n=Y_k+Y_{k+1}+\cdots+Y_n$. We have
%
%
\begin{eqnarray}
E|S_{2n}|^p&=&E|S_n+S_{n+1}^{2n}|^p\nonumber\\
&=&E\Biggl|\sum_{k=1}^n\bigl(Y_k-Y_k^{(n-k)}\bigr)+\sum_{k=1}^n
\bigl(Y_{n+k}-Y_{n+k}^{(k-1)}\bigr)\nonumber\\
&&\hspace*{70.2pt}{}
+\sum_{k=1}^nY_k^{(n-k)}+\sum_{k=1}^nY_{n+k}^{(k-1)}
\Biggr|^p\nonumber\\
&\leq&\Biggl(\sum_{k=1}^n\bigl\|Y_k-Y_k^{(n-k)}\bigr\|_p
+\sum_{k=1}^n\bigl\|Y_{n+k}-Y_{n+k}^{(k-1)}\bigr\|_p \\
&&\hspace*{66.4pt}{}+\Biggl\|\sum_{k=1}^nY_k^{(n-k)}
+\sum_{k=1}^nY_{n+k}^{(k-1)}\Biggr\|_p\Biggr)^p
\nonumber\\
&\leq&\Biggl(2D_p+\Biggl\|\sum_{k=1}^nY_k^{(n-k)}
+\sum_{k=1}^nY_{n+k}^{(k-1)}\Biggr\|_p\Biggr)^p
\nonumber\\
\label{eq:1_in_prop}
&=&\!:(2D_p+\|Z_n+W_n\|_p)^p.
\end{eqnarray}
For some positive constants $\psi_p$ that will be specified later, we
choose $C_p$ so that $C_p^{1/p}>D_p/\psi_p$. Then if $n\leq n_0$
\begin{eqnarray*}
E|Z_{n}|^p&\leq&(\|S_n\|_p+\|S_n-Z_n\|_p)^p\\
&\leq&(\|S_n\|_p+D_p)^p\\
&\leq&(1+\psi_p)^pC_p n^{p/2}.
\end{eqnarray*}
By the induction assumption, this relation holds with arbitrary $n$ for all
integer moments of order $\leq p-1$.
The same estimate applies for $E|W_n|^p$. Due to assumption~(B)
in Definition \ref{d}, the random variables
$Z_n$ and $W_n$ are independent. Thus
%
%
\begin{eqnarray}\label{eq:2_in_prop}
&&E|Z_n+W_n|^p\nonumber\\[-2pt]
&&\qquad\leq E|Z_n|^p+E|W_n|^p
+\sum_{m=1}^{p-1}\pmatrix{p\cr m}
E|Z_n|^mE|W_n|^{p-m}\nonumber\\[-10pt]\\[-10pt]
&&\qquad\leq n^{p/2}\Biggl[2(1+\psi_p)^pC_p
+\sum_{m=1}^{p-1}\pmatrix{p\cr m}
(1+\psi_m)^m(1+\psi_{p-m})^{p-m}C_mC_{p-m}\Biggr]\nonumber\\
&&\qquad=:n^{p/2}[2(1+\psi_p)^pC_p+R_p].\nonumber
\end{eqnarray}
Hence (\ref{eq:1_in_prop}) and (\ref{eq:2_in_prop}) and our
assumptions on $C_p$ imply that
%
%
\begin{eqnarray}\label{eq:3_in_prop}
E|S_{2n}|^p&\leq&\bigl(2\psi_pC_p^{1/p}+n^{1/2}[2(1+\psi_p)^pC_p
+R_p]^{1/p}\bigr)^p\nonumber\\[-10pt]\\[-10pt]
&\leq& C_p
n^{p/2}\bigl(2\psi_p+[2(1+\psi_p)^p+R_p/C_p]^{1/p}\bigr)^p.\nonumber
\end{eqnarray}
Choosing $\psi_p$ small enough, and then choosing $C_p$ large enough,
we can always achieve that the term in brackets of
(\ref{eq:3_in_prop}) is $\le\sqrt{2}$, provided that $p>2$, and that
the inequality $C_p^{1/p} > D_p/\psi_p$ mentioned before is satisfied.
Hence we have for every $n\leq n_0$
that $E|S_{2n}|^p\leq C_p(2n)^{p/2}$, proving (\ref{eq:ri0}) for all
even numbers $n\le2n_0$. The case of odd $n$ is similar. The proof of
Proposition~\ref{pr:mi} is finished for integer $p$.

For general $p>2$ we have by the result shown before
that (\ref{eq:ri}) holds for~$\lfloor p\rfloor$. (As usual,
$\lfloor~p\rfloor$ denotes the integer part of the real
number~$p$.) To finish the proof we need the following inequality which
will be
proven below:
%
%
\begin{eqnarray}\label{binineq}
|a+b|^p&\leq&|a|^p+|b|^p\nonumber\\[-10pt]\\[-10pt]
&&{}+\sum_{k=1}^{\lfloor p\rfloor}\pmatrix{p\cr k}
(|a|^k|b|^{p-k}+|b|^k|a|^{p-k}),\qquad
p\in[1,\infty).\nonumber\vadjust{\goodbreak}
\end{eqnarray}
Using (\ref{binineq}) we get a similar estimate for $E|Z_n+W_n|^p$ as in
(\ref{eq:2_in_prop}) and the proof can be finished along
the same lines as for integer $p$.

\textit{Verification of} (\ref{binineq}):
Let $x\in[0,1]$. We recall that
$(1+x)^p$ can be expanded in the binomial series
\[
(1+x)^p=\sum_{k\geq0}\pmatrix{p\cr k}x^k
\]
with
%
%
\begin{equation}\label{binomial}
\pmatrix{p\cr k}=\frac{p(p-1)\cdots(p-k+1)}{k!} .
\end{equation}
From (\ref{binomial}) it is clear that for $k\geq\lfloor p\rfloor+2$
we have
$\operatorname{sign}\{{p\choose k}\}=(-1)^{k-\lfloor p\rfloor+1}$.
This immediately yields for $k=\lfloor p\rfloor+2\ell$ with $\ell
\geq1$,
\[
\pmatrix{p\cr k}x^k+\pmatrix{p\cr k+1}x^{k+1}\leq
\pmatrix{p\cr k}x^{k}+\pmatrix{p\cr k+1}x^{k}
=
\pmatrix{p+1\cr k+1} x^{k}<0.
\]
Consequently
\[
\sum_{k\geq\lfloor p\rfloor+2}\pmatrix{p\cr k}x^k<0
\]
and
%
%
\begin{equation}\label{binser}
(1+x)^p\leq\sum_{k=0}^{\lfloor p\rfloor+1}\pmatrix{p\cr k}x^k.
\end{equation}
Now consider $|a+b|^p$.
If $|a|\geq|b|$, then we infer from (\ref{binser})
that
\begin{eqnarray*}
|a+b|^p&\leq&|a|^p\biggl(1+\biggl|\frac{b}{a}\biggr|\biggr)^p
\leq|a|^p \sum_{k=0}^{\lfloor p\rfloor+1}\pmatrix{p\cr k}\biggl|\frac
{b}{a}\biggr|^k\\
&=& |a|^p +\sum_{k=1}^{\lfloor p\rfloor}\pmatrix{p\cr k}|b|^k|a|^{p-k}\\
&&{}+ \pmatrix{p\cr\lfloor p\rfloor+1}|b|^p\biggl|\frac{b}{a}\biggr|^{\lfloor
p\rfloor+1-p}.
\end{eqnarray*}
Thus (\ref{binineq}) follows from ${p\choose\lfloor p\rfloor+1}
|\frac{b}{a}|^{\lfloor p\rfloor+1-p}\leq1$. Interchanging the
roles of $a$ and $b$ completes the proof.
\end{pf}

Using M\'{o}ricz \cite{moricz}, Theorem 1, we get:
\begin{Corollary}\label{co:mi}
Under the assumptions of Proposition \ref{pr:mi} with $p>2$, we have
for any $2<q\leq p$ and any $n\in\mathbb N$, $b\in\mathbb{Z}$
\[
E\max_{1\leq k\leq n}\Biggl|\sum_{j=b+1}^{b+k}Y_j\Biggr|^q\leq
C'_{p,q}
n^{q/2},
\]
where the constants $C'_{p,q}$ only depend on $p, q$ and the sequence
$\{Y_k\}$.
\end{Corollary}

A slightly weaker result can also be derived from
Proposition \ref{pr:mi} for the case of $0<q\le2$.

\subsection{\texorpdfstring{Proofs of Theorems \protect\ref{th:poly} and \protect\ref{th:exp}}
{Proofs of Theorems 1 and 2}}

We give the proof of
Theorem \ref{th:poly}.
Note first of all that $\delta(m)=
\|Y_k-Y_k^{(m)}\|_p\geq\|Y_k-Y_k^{(m)}\|_2$, and consequently
(\ref{eq:sighelp}) holds when the $L^2$-norm is replaced by the $L^p$-norm.
Since $A>1$ in (\ref{Acond}), we infer that the series in
(\ref{eq:sig}) is absolutely convergent.

Let us specify some constants that will be used for the proof. By
our assumption on $A$ it is possible to find a constant
$0<\varepsilon_0<1/2$ such that
\[
A>\frac{p-2}{2\eta(1-\varepsilon_0)^2}\biggl(1-\frac{1+\eta
}{p}\biggr).
\]
Then we set
%
%
\begin{eqnarray}\label{eq:dab}
\delta&=&\frac{\beta}{1+\alpha} \quad\mbox{with}\quad \alpha=
\frac{2\eta(1-\varepsilon_0)}{p-2(1+\eta)},\nonumber\\[-8pt]\\[-8pt]
\beta&=&(1-\varepsilon_0)\alpha.\nonumber
\end{eqnarray}
For some $\varepsilon_1>0$ (which will be specified later)
we now define $m_k=\lfloor
\varepsilon_1k^\delta\rfloor$. The first step in the proof of
(\ref{eq:strappr}) is to show that it is sufficient to provide the
strong approximation for the perturbed sequence $Y_k'=Y_k^{(m_k)} $. We
notice that our main assumption (\ref{eq:ma}) yields
$\|Y_{k}-Y_k'\|_p\ll k^{-A\delta}$. If $A\delta<1$, then
\begin{eqnarray*}
&&P\Biggl(\max_{2^n\leq k\leq
2^{n+1}}\Biggl|\sum_{j=1}^k(Y_j-Y_j')\Biggr|>
\frac{1}{n}2^{({n}/{p})(1+\eta)}\Biggr)\\
&&\qquad \leq P\Biggl(\sum_{j=1}^{2^{n+1}}|Y_j-Y_j'|>
\frac{1}{n}2^{({n}/{p})(1+\eta)}\Biggr)\\
&&\qquad \leq2^{-n(1+\eta)}n^p\Biggl(\sum_{j=1}^{2^{n+1}}\|Y_j-Y_j'\|
_p\Biggr)^p\\
&&\qquad \ll2^{-c_1n}n^p,
\end{eqnarray*}
where $c_1=(1+\eta)-(1-A\delta)p>0$. Thus by the Borel--Cantelli
lemma we have almost surely
\[
\sum_{j=1}^kY_j=\sum_{j=1}^kY_j'+o\bigl(k^{(1+\eta)/p}\bigr)
\qquad\mbox{a.s.}
\]
If $A\delta\geq1$ we get an (even better) error term of order
$o(k^{1/p})$.

The main part of the proof of Theorem \ref{th:poly} is based on a
blocking argument. We partition $\mathbb{N}$ into disjoint blocks
\[
\mathbb{N}=J_1\cup I_1\cup J_2\cup I_2\cup\cdots,
\]
where $|I_k|=\lfloor k^{\alpha}\rfloor$ and $|J_k|=\lfloor
k^{\beta}\rfloor$ with $\alpha$, $\beta$ as in (\ref{eq:dab}). Let
us further set
\[
I_k=\{\underlineIk,\ldots,\overline{i}_k\} \quad\mbox{and}\quad
J_k=\{\underlineJk,\ldots,\overline{j}_k\}
\]
and
\[
\xi_k=\sum_{j\in I_k}Y_j' \quad\mbox{and}\quad \eta_k=\sum_{j\in
J_k}Y_j'.
\]
Note that $\overline{i}_k=O(k^{1+\alpha})$. Provided
that $\varepsilon_1$ in the definition of $m_k$ is chosen small
enough, this will imply that
\[
|J_k|=\lfloor k^{\beta}\rfloor> \lfloor
\varepsilon_1\underlineIk^{\delta} \rfloor= m_{\underlineIk},
\]
and hence by assumption (B) it follows that $\{\xi_k\}$ and $\{\eta
_k\}$
each define a~sequence of independent random variables.

The following lemma by Sakhanenko \cite{sakhanenko} (cf. also Shao
\cite{shao}) is our crucial ingredient for the
construction of the approximating processes.
\begin{Lemma}\label{le:sak}
Let $\{\xi_k\}$ be a sequence of centered independent random
variables with finite $p$th moments, $p>2$. Then we can redefine
$\{\xi_k\}$ on a suitable probability space, together with a
sequence $\{\xi_k^*\}$ of independent normal random variables with
$E\xi_k^*=0$, $E(\xi_k^*)^2=E\xi_{k}^2$ such that for any $x>0$,
$m\geq1$
\[
P\Biggl(\max_{1\leq k\leq
m}\Biggl|\sum_{j=1}^k\xi_j-\sum_{j=1}^k\xi_j^*\Biggr|>x
\Biggr)\leq
C\frac{1}{x^p}\sum_{j=1}^mE|\xi_j|^p,
\]
where $C$ is an absolute constant.
\end{Lemma}

We shall now apply Lemma \ref{le:sak} to the sequences $\{\xi_k\}$
and $\{\eta_k\}$. For this purpose we need estimates of the moments
$E|\xi_k|^p$, $E|\eta_k|^p$. By Minkowski's inequality and
Proposition \ref{pr:mi} we get
\begin{eqnarray*}
E|\xi_k|^p&\leq&
\biggl(\biggl\|\sum_{j\in I_k}Y_k\biggr\|_p+
\sum_{j\in I_k}\|Y_j-Y_j'\|_p\biggr)^p\\
&=&O\bigl((|I_k|^{1/2}+|I_k|\cdot\underlineIk^{-A\delta
}
)^p\bigr).
\end{eqnarray*}
Some easy algebra shows that the restrictions on the parameters $A$,
$\delta$, $\alpha$ and~$\varepsilon_0$ imply
\[
|I_k|\cdot\underlineIk^{-A\delta}\ll k^{\alpha}\cdot
k^{-A\delta(1+\alpha)}\ll k^{\alpha/2}\ll|I_k|^{1/2}.
\]
A similar estimate holds for $E|\eta_k|^p$. Hence we can find
constants $F_p$ such that
\[
E|\xi_k|^p\leq F_p|I_k|^{p/2}\vadjust{\goodbreak}
\]
and
\[
E|\eta_k|^p\leq
F_p|J_k|^{p/2},
\]
where $F_p$ does not depend on $k$.

Let $L_n=\sum_{k=1}^{n}|I_k|$. Then $L_n=O(n^{(1+\alpha)})$.
By our previous estimates and by Lemma \ref{le:sak} we infer that,
after enlarging the probability space, we have
%
%
\begin{eqnarray}\label{eq:sakh_poly}
&&P\Biggl(\max_{2^n\leq k\leq
2^{n+1}}\Biggl|\sum_{j=1}^k\xi_j-\sum_{j=1}^k\xi_j^*
\Biggr|>L_{2^n}^{({1+\eta})/{p}}\Biggr) \nonumber\\
&&\qquad\leq L_{2^n}^{-(1+\eta
)}\sum
_{k=1}^{2^{n+1}}E|\xi_k|^p\\
&&\qquad=O\bigl(2^{[-(1+\alpha)(1+\eta)+{\alpha
p}/{2}+1]n}\bigr),\nonumber
\end{eqnarray}
where $\xi_k^*$ is a sequence of independent and centered normal
random variables with $E(\xi_k^*)^2=E\xi_k^2$.
The
exponent in (\ref{eq:sakh_poly}) will be negative if
$(1+\alpha)(1+\eta)>\frac{\alpha p}{2}+1$. This is equivalent to
$\alpha<\frac{2\eta}{p-2(1+\eta)}$, which follows by (\ref{eq:dab}).
Thus, by the Borel--Cantelli lemma we obtain
\[
\sum_{j=1}^k\xi_j=\sum_{j=1}^k\xi_j^*
+O\bigl(L_k^{({1+\eta})/{p}}\bigr) \qquad\mbox{a.s.}
\]
By further enlarging the probability space we can write
\[
\sum_{j=1}^k\xi_j=W_1\Biggl(\sum_{j=1}^k\operatorname{Var}(\xi_j)\Biggr)
+O\bigl(L_k^{({1+\eta})/{p}}\bigr) \qquad\mbox{a.s.},
\]
where $\{W_1(t), t\geq0\}$ is a standard Wiener process. The same
arguments show that
\[
\sum_{j=1}^k\eta_j=W_2\Biggl(\sum_{j=1}^k\operatorname{Var}(\eta
_j)\Biggr)
+O\bigl(M_k^{({1+\eta})/{p}}\bigr) \qquad\mbox{a.s.},
\]
where $\{W_2(t), t\geq0\}$ is another standard Wiener process on
the same probability space and $M_n={\sum_{k=1}^{n}}|J_k|$.

We define
\[
b_k^2=\operatorname{Var}\biggl(\sum_{j\in
I_k}Y_j'\biggr)\Big/|I_k|
\]
and
\[
h_k^2=\operatorname{Var}\biggl(\sum_{j\in J_k}Y_j'\biggr)\Big/|J_k|.
\]
For $\ell\in I_k$ we
set $\sigma_{\ell}^2=b_k^2$ and for $\ell\in J_k$ we set $\sigma
_{\ell}^2=0$.
Similarly define $\tau_{\ell}^2=h_k^2$ if $\ell\in J_k$ and $\tau
_{\ell
}^2=0$ if $\ell\in I_k$. Put
\[
s_n^2=\sum_{k=1}^n \sigma_k^2,\qquad t_n^2=\sum_{k=1}^n \tau_k^2\qquad
(n=1, 2, \ldots).
\]
Summarizing our results so far we can write
\[
\sum_{k=1}^{\overline{i}_n}Y_k=W_1\Biggl(\sum_{k=1}^{\overline
{i}_n}\sigma_k^2\Biggr)+
W_2\Biggl(\sum_{k=1}^{\overline{i}_n}\tau_k^2\Biggr)+O
\bigl(\overline
{i}_n^{(1+\eta)/p}\bigr) \qquad\mbox{a.s.}
\]
It is a basic result that our stationarity and dependence assumptions\break imply
%
%
\begin{eqnarray}\label{eq:assig}
\operatorname{Var}\biggl(\sum_{j\in
I_k}Y_j\biggr)\Big/|I_k|&=&\sigma^2+O(k^{-\xi})
\quad\mbox{and}\nonumber\\[-8pt]\\[-8pt]
\operatorname{Var}\biggl(\sum_{j\in J_k}Y_j\biggr)\Big/|J_k|&=&\sigma
^2+O(k^{-\xi})\nonumber
\end{eqnarray}
as $k\to\infty$, for some small enough $\xi>0$.
It can be easily shown that (\ref{eq:assig}) remains true if
the $Y_j$ are replaced with $Y_j'$.
Indeed, by the Minkowski inequality we infer that
\begin{eqnarray*}
\operatorname{Var}^{1/2}\biggl(\sum_{j\in I_k}Y_j'\biggr)&\leq&
\operatorname{Var}^{1/2}\biggl(\sum_{j\in I_k}Y_j\biggr)+
\operatorname{Var}^{1/2}\biggl(\sum_{j\in I_k}(Y_j-Y_j')\biggr)\\
&\leq&\operatorname{Var}^{1/2}\biggl(\sum_{j\in I_k}Y_j\biggr)+|I_k|\max_{j\in
I_k}\|Y_j-Y_j'\|_2.
\end{eqnarray*}
Furthermore, using the definitions of the introduced constants we obtain
\begin{eqnarray*}
\max_{j\in I_k}\|Y_j-Y_j'\|_2
&\ll&\underlineIk^{-A\delta}
\ll k^{-(\alpha+1)A\delta}\\
&\leq& k^{-\beta}=k^{-(1-\epsilon_0)\alpha}\qquad \mbox{with
$\epsilon_0<1/2$}.
\end{eqnarray*}
Since by definition $|I_k|\ll k^{\alpha}$, we conclude that
\[
\operatorname{Var}^{1/2}\biggl(\sum_{j\in I_k}Y_j'\biggr)\Big/|I_k|^{1/2}\leq
\operatorname{Var}^{1/2}\biggl(\sum_{j\in I_k}Y_j
\biggr)\Big/|I_k|^{1/2}+O
\bigl(k^{\alpha(\epsilon_0-1/2)}\bigr)
\]
as $k\to\infty$.
In the same manner a lower bound for $\operatorname{Var}^{1/2}(\sum
_{j\in I_k}Y_j')/|I_k|^{1/2}$
can be obtained. Proving the analogue of the second part of (\ref
{eq:assig}) for the $Y_j'$ is similar.

In other words, we have shown (\ref{eq:strappr}) along the subsequence
$\{\overline{i}_n\}$
with values of $s_n^2$ and $t_n^2$ that satisfy (\ref{eq:sn})
and (\ref{stlimsup}).
The relation $|s_n^2-\sigma^2n|=O(n^{1-\epsilon})$, $\epsilon>0$,
follows by simple calculations.

To finish the proof we have to show that the fluctuations of the
partial sums and the Wiener processes $W_1$ and $W_2$
within the blocks $I_k$ are small enough. Since fluctuation properties of
Wiener processes are easy to handle using standard deviation
inequalities (see, e.g., \cite{cr}), we only investigate
the partial sums. By
Corollary~\ref{co:mi} we have
\begin{eqnarray*}
P\Biggl(\sup_{\underlineIk\leq\ell\leq
\overline{i}_k}\Biggl|\sum_{j=\underlineIk}^{\ell}
Y_j\Biggr|>\underlineIk^{({1+\eta})/{p}}\Biggr)
&\leq&\underlineIk^{-(1+\eta)} E\Biggl(\sup_{\underlineIk\leq
\ell\leq
\overline{i}_k}\Biggl|\sum_{j=\underlineIk}^{\ell}Y_j
\Biggr|^p\Biggr)\\
&\ll&\underlineIk^{-(1+\eta)}|I_k|^{p/2}\\
&\ll& k^{-(1+\eta)(1+\alpha)+{\alpha p}/{2}}\\
&=&O\bigl(k^{-(1+\varepsilon_2)}\bigr),
\end{eqnarray*}
if $\varepsilon_2>0$ is chosen sufficiently small. The
Borel--Cantelli lemma shows that we can also control the fluctuation
within the blocks. Thus (\ref{eq:strappr}) is proven.

The proof of Theorem \ref{th:exp} is similar to the proof of
Theorem \ref{th:poly} and will be therefore omitted.
We only remark that under the exponential mixing
rate logarithmic block sizes are required in the blocking
argument.



\subsection{\texorpdfstring{Proof of Proposition \protect\ref{pr:corr}}
{Proof of Proposition 1}}

We use the notation introduced in the proof of Theorem \ref{th:poly}.
Further we let $I=I_1\cup I_2 \cup\cdots$ and
$J=J_1\cup J_2\cup\cdots$ and $M_n=\{1,\ldots,n\}$. By
looking at the proof of Theorem \ref{th:poly}, it readily follows
that
\begin{eqnarray*}
\frac{1}{s_n}\sum_{i\in I\cap M_n}Y_i&=&W_1(s_n^2)/s_n-X_n,\\
\frac{1}{t_m}\sum_{j\in J\cap M_n}Y_j&=&W_2(t_m^2)/t_m-Z_m,
\end{eqnarray*}
where
%
%
\begin{eqnarray}\label{error}
X_n&=&o\bigl((s_n^2)^{({1+\eta})/{p}-{1}/{2}}\bigr)=o(1)
\qquad\mbox
{a.s.} \quad\mbox{and}\nonumber\\[-8pt]\\[-8pt]
Z_m&=&o\bigl((t_m^2)^{({1+\eta})/{p}-{1}/{2}}\bigr)=o(1)
\qquad\mbox{a.s.}\nonumber
\end{eqnarray}
Hence
\begin{eqnarray*}
&&\corr(W_1(s_n^2),W_2(t_m^2))\\
&&\qquad =\corr\biggl(\frac{1}{s_n}W_1(s_n^2),
\frac{1}{t_m}W_2(t_m^2)\biggr)\\
&&\qquad =\corr\biggl(\frac{1}{s_n}\sum_{i\in I\cap M_n}Y_j+X_n,
\frac{1}{t_m}\sum_{j\in J\cap M_n}Y_j+Z_m\biggr).
\end{eqnarray*}
In order to calculate this correlation we need
a couple of estimates.

First we note that by the definition of $s_n^2$ and $t_n^2$
%
%
\begin{equation}
s_n^2\sim\sigma^2|I\cap M_n| \quad\mbox{and}\quad
t_m^2\sim\sigma^2|J\cap M_m|.
\end{equation}
It readily follows from Proposition \ref{pr:mi} that
%
%
\begin{equation}\label{eq:permom}
\biggl\|\frac{1}{s_n}\sum_{i\in I\cap M_n}Y_i\biggr\|_p\leq C_p,
\end{equation}
where $C_p$ does not depend on $n$.
Thus
\begin{eqnarray*}
\sup_{n\geq1}\|X_n\|_p&=&\sup_{n\geq1}\biggl\|
\frac{1}{s_n}\sum_{i\in I\cap M_n}Y_i-W_1(s_n^2)/s_n\biggr\|_p\\
&\leq&\sup_{n\geq1}\biggl\|
\frac{1}{s_n}\sum_{i\in I\cap M_n}Y_i\biggr\|_p+\|W_1(1)\|
_p<\infty,
\end{eqnarray*}
and hence $\{X_n^2\}$ is uniformly integrable. This and (\ref{error})
show that
$\Var(X_n)\to0$; by the same arguments $\Var(Z_m)\to0$.
By (\ref{eq:assig})
%
%
\begin{equation}\label{eq:permom2}
\biggl\|\frac{1}{s_n}\sum_{i\in I\cap M_n}Y_i\biggr\|_2\sim\sigma^2.
\end{equation}
Thus by (\ref{eq:permom}) and (\ref{eq:permom2})
\begin{eqnarray*}
c_{1}(m,n):\!&=&\Cov\biggl(
Z_m,\frac{1}{s_n}\sum_{i\in I\cap M_n}Y_i
\biggr)\\
&\leq&\Var^{1/2}(Z_m)\Var^{1/2}\biggl(\frac{1}{s_n}\sum_{i\in I\cap M_m}Y_i
\biggr)\\
&=&o(1) \qquad\mbox{for $m,n\to\infty$},
\end{eqnarray*}
and similarly
\[
c_2(m,n):=\Cov\biggl(
X_n,\frac{1}{t_m}\sum_{j\in J\cap M_m}Y_j
\biggr)=o(1) \qquad\mbox{for $m,n\to\infty$.}
\]
Furthermore we have
\begin{eqnarray*}
B_{1}(n):\!&=&\Var^{1/2}\biggl(\frac{1}{s_n}\sum_{i\in I\cap M_n}Y_i+X_n
\biggr)\\
&\geq&\Var^{1/2}\biggl(\frac{1}{s_n}\sum_{i\in I\cap M_n}Y_i
\biggr)-\Var^{1/2}(X_n)\\
&=&\sigma+o(1) \qquad\mbox{for $n\to\infty$}
\end{eqnarray*}
and
\begin{eqnarray*}
B_2(m):\!&=&\Var^{1/2}\biggl(\frac{1}{t_m}\sum_{j\in J\cap M_m}Y_j+Z_m
\biggr)\\
&\geq&\sigma+o(1) \qquad\mbox{for $m\to\infty$}.
\end{eqnarray*}
Finally we introduce the term
\[
c_0(m,n)=\frac{1}{s_nt_m}\sum_{i\in I\cap M_n}\sum_{j\in J\cap
M_m}\Cov
(Y_i,Y_j).
\]
We choose $r\geq0$ such that $n\in I_{r+1}\cup J_{r+1}$, and we choose
$v\geq0$ such that $m\in I_{v+1}\cup J_{v+1}$ and recall that
by Theorem \ref{th:poly} we have
${\sum_{i\in\mathbb{Z}}}|{\Cov(Y_0,Y_i)}|<\infty.
$
Hence if $v\leq2r$ we have
\begin{eqnarray*}
c_0(m,n)
&\leq& {s_{\overline{i}_{r}}^{-1}t_{\overline{i}_{v}}^{-1}\sum_{i\in
I_1\cup\cdots\cup I_{r+1}}\sum_{j\in J_1\cup\cdots\cup
J_{v+1}}}|\Cov
(Y_i,Y_j)|\\
&\leq& s_{\overline{i}_{r}}^{-1}t_{\overline{i}_{v}}^{-1}\sum_{j\in
J_1\cup\cdots\cup J_{v+1}}\sum_{i\in\mathbb{Z}}|\Cov
(Y_i,Y_j
)|\\
&=& s_{\overline{i}_{r}}^{-1}t_{\overline{i}_{v}}^{-1}\sum_{j\in
J_1\cup
\cdots\cup J_{v+1}}\sum_{i\in\mathbb{Z}}|\Cov(Y_i,Y_0)|\\
&\ll& s_{\overline{i}_{r}}^{-1}t_{\overline{i}_{v}}^{-1}
(|J_1|+\cdots
+ |J_{v+1}|)\\
&\ll& s_{\overline{i}_{r}}^{-1}t_{\overline{i}_{v}}^{-1}t_{\overline
{i}_{v+1}}^{2}
=o(1) \qquad\mbox{as $m,n\to\infty$.}
\end{eqnarray*}
If $v> 2r$, we have to additionally show that
\[
{s_{\overline{i}_{r}}^{-1}t_{\overline{i}_{v}}^{-1}\sum_{i\in I_1\cup
\cdots\cup I_{r+1}}\sum_{j\in J_{2r+1}\cup\cdots\cup J_{v+1}}}|\Cov
(Y_i,Y_j)|\to0.
\]
Now we have by (\ref{eq:sighelp}) and assumptions (\ref{eq:ma}),
(\ref{Acond})
that
\begin{eqnarray*}
&&{s_{\overline{i}_{r}}^{-1}t_{\overline{i}_{v}}^{-1}\sum_{i\in
I_1\cup
\cdots\cup I_{r+1}}\sum_{j\in J_{2r+1}\cup\cdots\cup J_{v+1}}}|\Cov
(Y_i,Y_j)|\\
&&\qquad \leq s_{\overline{i}_{r}}^{-1}t_{\overline{i}_{v}}^{-1}\sum
_{\pi
\geq2r+1}
|J_{\pi}|\sum_{\ell=1}^{r+1}|I_\ell|(d(I_\ell,J_{\pi})
)^{-1}\\
&&\qquad \ll s_{\overline{i}_{r}}^{-1}t_{\overline{i}_{v}}^{-1}\sum
_{\pi
\geq2r+1}
\pi^{\beta}\sum_{\ell=1}^{r+1}\ell^{\alpha}(d(I_\ell,J_{\pi
}))^{-1}.
\end{eqnarray*}
For $\ell\in\{1,\ldots,r+1\}$ and $\pi\geq2r+1$ we have constants
$k_0$ and $k_1$ independent of $r$ and $\pi$ such that
\[
d(I_\ell, J_{\pi})\geq k_0( \pi^{\alpha+1}-r^{\alpha+1}
)\geq
k_1\pi^{\alpha+1}
\]
and thus
\begin{eqnarray*}
&&s_{\overline{i}_{r}}^{-1}t_{\overline{i}_{v}}^{-1}\sum_{\pi\geq2r+1}
\pi^{\beta}\sum_{\ell=1}^{r+1}\ell^{\alpha}(d(I_\ell,J_{\pi
})
)^{-1}\\
&&\qquad \ll s_{\overline{i}_{r}}^{-1}t_{\overline
{i}_{v}}^{-1}r^{\alpha
+1}\sum_{\pi\geq2r+1}
\pi^{\beta-\alpha-1} \ll s_{\overline{i}_{r}}^{-1}t_{\overline{i}_{v}}^{-1}
r^{1+\beta}\\
&&\qquad\ll r^{-(\alpha-\beta)/2}=o(1) \qquad\mbox{as $r\to\infty$}.
\end{eqnarray*}
Using the definitions of $c_0$, $c_1$, $c_2$ and $B_1$ and $B_2$ we
see that
%
%
\begin{equation}\label{cov}
\corr(W_1(s_n^2),W_2(t_m^2))
=\frac{
c_0(m,n)+c_1(m,n)+c_2(m,n)+\operatorname{cov}(X_n,Z_m)}{B_1(n)B_2(m)}.\hspace*{-28pt}
\end{equation}
We have shown
$c_0(m,n)+c_1(m,n)+c_2(m,n)+\operatorname{cov}(X_n,Z_m)\to0$ as
$m,n\to\infty$ while
the denominator in (\ref{cov}) is bounded away from zero. This finishes
the proof of Proposition \ref{pr:corr}.

\section*{Acknowledgments}
The authors are indebted to Wei Biao Wu for raising the question
leading to the present research and for valuable comments. We also
thank two anonymous referees for
several inspiring questions and remarks which led to a considerable
improvement of the presentation.

%

%
\printaddresses

\end{document}